\documentclass{amsart}
\usepackage[dvips]{graphics}

\theoremstyle{definition}

\theoremstyle{remark}

\numberwithin{equation}{section}

\theoremstyle{problem}

\newcommand{\sgn}{\mathop{\rm sgn}\nolimits}

\newcommand{\R}{\mathop{\rm Re}\nolimits}
\newcommand{\const}{\mbox{const}}

\newcommand{\Ga}{\alpha}
\newcommand{\Gb}{\beta}
\newcommand{\Gd}{\delta}

\newcommand{\Gve}{\varepsilon}
\newcommand{\Gf}{\phi}
\newcommand{\Gvf}{\varphi}
\newcommand{\Gg}{\gamma}
\newcommand{\Gc}{\chi}

\newcommand{\Gk}{\kappa}

\newcommand{\Gl}{\lambda}
\newcommand{\Gn}{\eta}

\newcommand{\Gs}{\sigma}

\newcommand{\Gx}{\xi}

\newcommand{\Gz}{\zeta}

\newcommand{\GF}{\Phi}
\newcommand{\GG}{\Gamma}

\newcommand{\beq}{\begin{equation}}
\newcommand{\eeq}{\end{equation}}
\newcommand{\barr}{\begin{eqnarray}}
\newcommand{\earr}{\end{eqnarray}}
\newcommand{\beqn}{\begin{equation*}}
\newcommand{\eeqn}{\end{equation*}}
\newcommand{\barrn}{\begin{eqnarray*}}
\newcommand{\earrn}{\end{eqnarray*}}

\newcommand{\fr}{\frac}

\begin{document}

\title[Semi-infinite Hilbert transform and applications] {Integral relations associated with the semi-infinite Hilbert transform and applications to singular integral equations}
\author{Y.A. Antipov}
\address{Department of Mathematics, Louisiana State University, Baton
Rouge, Louisiana 70803}
\email{yantipov@lsu.edu}
\thanks{Research of the first author was sponsored by the Army Research Office and was accomplished under Grant Number {\bf W911NF-17-1-0157}.
The views and conclusions contained in this document are those of the authors and should not be interpreted 
as representing the official policies, either expressed or implied, of the Army Research Office or the U.S. Government.
The U.S. Government is authorized to reproduce and distribute reprints for Government purposes notwithstanding any copyright notation herein.}

\author{S.M. Mkhitaryan}
\address{Department of Mechanics of Elastic and Viscoelastic Bodies, National Academy of Sciences,
Yerevan 0019, Armenia}
\email{smkhitaryan39@rambler.ru}


\begin{abstract}

Integral relations with the Cauchy kernel on a semi-axis for the Laguerre polynomials, the confluent hypergeometric function, and
the cylindrical functions are derived. A part of these formulas are obtained by exploiting some properties
of the Hermite polynomials including their Hilbert and Fourier transforms and connections to the Laguerre polynomials.   
The relations discovered give rise to complete systems of new orthogonal functions. Free of singular integrals exact and approximate solutions  to
 the characteristic and complete singular integral
equations in a semi-infinite interval are proposed. Another set of the Hilbert transforms in a semi-axis are deduced from  integral relations with the Cauchy kernel in a finite
segment for the Jacobi polynomials and the Jacobi functions of the second kind by letting some parameters involved go to infinity.
These formulas lead to integral relations for the Bessel functions. Their application to a model problem of contact mechanics is given.
A new quadrature formula for the Cauchy integral in a semi-axis based on an integral relation for the Laguerre polynomials and 
the confluent hypergeometric function is derived and tested numerically. Bounds for the reminder are found.

\end{abstract}

\maketitle

\setcounter{equation}{0}

\section{Introduction} 

The Hilbert transform, convolution of a function and the Cauchy kernel, in the real axis was introduced by Hilbert in the beginning
of the twentieth  century as a tool for boundary value problems of the theory of analytic functions. Since then its properties and related methods
in miscellaneous areas of applied sciences have been of interest for pure and applied mathematicians.
Properties of the Hilbert transforms in the real axis are presented in detail in \cite{tit}. The Hilbert
transform in a finite segment has been extensively studied due to applications to singular integral equations 
 \cite{rei},  \cite{ell},  \cite{mus1},  \cite{son},  \cite{tri2},  \cite{gak},
 and quadrature formulas \cite{kor}, \cite{sta} 
for the Cauchy integral in the interval $(-1,1)$. In particular, a formula for the Hilbert transform in the interval $(-1,1)$
of the weighted Jacobi polynomials was discovered in \cite{tri1}.  This formula generates many integral relations for orthogonal polynomials crucial  for the solution of  singular integral equations. In \cite{kop},  the Fourier-Plancherel  transformation was employed to derive a spectral representation of the finite Hilbert transform and an expansion of an arbitrary $L^2(a,b)$-function in terms of the Hilbert operator eigenfunctions.  In \cite{mkh1}, 
a self-adjoint differential equation was  considered and integral operators  whose eigenfunctions are the solutions
of that differential equation were determined. This made possible to recover spectral relations for the operators found and solve the corresponding singular integral equations in a finite interval.
Finite Hilbert transforms of the Chebyshev, Bernstein, and Lagrange interpolating polynomials were used for an approximate
solution of the Prandtl integro-differential equation in \cite{mkh2}.

The Hilbert transform in a semi-axis has received less attention than their infinite and finite analogues. At the same time, many model problems including 
those arising in fluid mechanics, fracture, and penetration mechanics are governed by singular integral equations or their systems in a semi-axis whose kernels can be represented as a sum
of the Cauchy kernel and a regular kernel.  Also,  vector Riemann-Hilbert problems, when the associated matrix Wiener-Hopf  
factors are infeasible, do not admit a closed-form solution. Alternatively, they
can be written as systems of singular integral equations.  The  kernels of these systems 
may often be split into the Cauchy kernel and bounded functions (see for example, \cite{ant1}, \cite{ant2}, \cite{ant3}).
Replacing the semi-axis $(0,\infty)$ by a finite segment $(0,A)$ and applying the numerical methods for singular integral equations in a finite segment
do not preserve the properties of the solution far away from the point $0$ and may generate a significant error of approximation due to the change of the weight function.
In this paper, we aim to derive a series of  integral relations  in a semi-axis for the Laguerre polynomials, the confluent hypergeometric function, the Tricomi function,
and the Bessel functions and study a new system of functions $\{G_m(x)\}_{m=0}^\infty$, the $G$-functions, generated by these integral relations.
We also intend, by means of these integral relations, to solve some singular integral equations and obtain a quadrature formula for the Cauchy integral in a semi-axis.

In Section 2, we introduce the $G$-functions, the Hilbert transforms of the weighted Hermite polynomials, 
\beq
G_n(x)=\fr{1}{\pi}\int_{-\infty}^\infty \fr{H_n(t)e^{-t^2/2}dt}{t-x}, \quad n=0,1,\ldots,
\label{1.1}
\eeq
show that they constitute a complete orthogonal system 
in the associated space, and derive their representation in terms of the confluent hypergeometric function $\GF$,
\beq
G_{2m+j}(x)=\fr{\sqrt{2}(2m+j)! x^{1-j}}{(-1)^{j+1} e^{x^2/2}}\sum_{k=0}^m\fr{ 2^k\GF(1/2-k-j,3/2-j;x^2/2)}{(-1)^k(m-k)!\GG(k+j+1/2)}, 
\label{1.2}
\eeq	
where  $m=0,1,\ldots$ and  $j=0,1$. We discover the Hilbert transforms in a semi-axis of the Laguerre polynomials 
in two weighted spaces, $L^2_{w_\pm}(0,\infty)$ and  $L^2_{v_\pm}(0,\infty)$, where $w_\pm(x)=x^{\pm 1/2} e^{-x/2}$ and $v_\pm(x)=x^{\pm 1/2} e^{-x}$.  It is found that the first group
of the relations generates an orthogonal basis whose elements are  the $G$-functions, while the second one gives rise to a nonorthogonal basis. 
We study the $G$-functions in Section 3. It is shown that these functions of even and odd indices satisfy certain second-order inhomogeneous differential equations         
with variable coefficients. Also, the asymptotics of the $G$-functions at the points $0$ and $\infty$  is derived. 

In Section 4,  based on the Tricomi  integral relation \cite{tri1} written for the Jacobi polynomials  $P_n^{(\Ga,\Gb)}(1-2x/\Gb)$ and by passing
to the limit $\Gb\to\infty$, we discover integral relations for the Laguerre polynomials $L_n^\Ga(x)$ 
to be used in Section 6 for a quadrature formula for the Cauchy integral in a semi-axis. In addition, by employing some representations
of the Jacobi function of the second kind $Q_n^{(\Ga,\Gb)}(x)$ and passing to the limit $n\to\infty$ in the representations for $n^{-\Ga}Q_n^{(\Ga,\Gb)}(1-\fr12z^2/n^2)$
we obtain integral relations for the Bessel functions. 

In Section 5, we apply the integral relations in a semi-axis derived in the previous sections  to find a closed-form solution of the integral equation with the Cauchy kernel in the interval $(0,\infty)$
in the classes of bounded and unbounded  at zero functions. In addition, we construct two bilinear expansions of the Cauchy kernel in terms of the Laguerre polynomials 
and the $G$-functions. We also solve a system of two complete integral equations with the Cauchy kernel in a semi-axis by 
reducing it  to
an infinite system  of linear algebraic equations of the second kind. By using the integral relation for the Bessel function obtained in  Section 4
and the Hankel transform   we obtain a closed-form solution, that is free of singular integrals, to a
contact problem on a semi-infinite stamp and an elastic half-plane.

 In Section 6, we obtain the following quadrature formula for the Cauchy integral  in a semi-axis:
\beq
 \fr{1}{\pi}\int_0^\infty \fr{f(t)t^\Ga e^{-t}dt}{t-x}=-\fr{1}{n+\Ga}\sum_{m=1}^n\fr{x_m  f(x_m)}{L_{n-1}^\Ga(x_m)}\fr{Q^\Ga_n(x)-Q^\Ga_n(x_m)}{x-x_m},\quad 0<x<\infty,
\label{1.3}
\eeq
where
\beq
Q^\Ga_n(x)=\fr1{\pi}e^{-x}\GG(\Ga)\GF(-n-\Ga,1-\Ga;x)-\cot\pi\Ga x^\Ga e^{-x}L_n^\Ga(x).
\label{1.4}
\eeq
It is exact for a polynomial of degree $n-1$ and requires $n$ zeros, $x_m$, of the Laguerre polynomial $L_n^\Ga(x)$.

\setcounter{equation}{0}

\section{Hilbert transforms of the weighted Hermite and Laguerre polynomials and the $G$- and $V$-functions}\label{form}

In this section we introduce the Hilbert transforms $G_n(x)$ and $V_n(x)$ of the weighted Hermite polynomials,  derive associated integral relations for the weighted Laguerre polynomials
and study the properties of the functions  $G_n(x)$ and $V_n(x)$.

\subsection{Relations for the weighted Laguerre polynomials $e^{-\Gn/2}\Gn^{\pm 1/2}L_m(\Gn)$}

For real functions $\Gvf_k(x)$ such that $\Gvf_k(x)\in L^2(-\infty,\infty)$, $k=1,2$, define their Hilbert transform $\GF_k(x)=H[\Gvf_k(t)](x)$ by
\beq
\GF_k(x)=\fr{1}{\pi}\int_{-\infty}^\infty\fr{\Gvf_k(t)dt}{t-x}, \quad -\infty<x<\infty.
\label{2.1}
\eeq
The integral is understood in the sense of the principal value at $t=x$ and the mean square sense at infinity. Then
$\GF_k(x)\in L^2(-\infty,\infty)$, $k=1,2$, and the generalized Parseval's relation holds  \cite{tit}
\beq
\int_{-\infty}^\infty\Gvf_1(x)\Gvf_2(x)dx=\int_{-\infty}^\infty\GF_1(x)\GF_2(x)dx.
\label{2.2}
\eeq
The Hilbert transform $H$ is a 1-1 map  and  a unitary operator in the space $L^2(-\infty,\infty)$. Since
this operator preserves the inner product in the Hilbert space $L^2(-\infty,\infty)$, a complete orthogonal system in the space $L^2(-\infty,\infty)$ is transformed
by the operator $H$ into another complete orthogonal system in the same space. 

\vspace{2mm}

{\sc Theorem 2.1} Let $H_n(x)$ be the Hermite polynomials normalized by the condition  
\beq
\lim_{x\to\infty}x^{-n}H_n(x)=2^n.
\label{2.2'}
\eeq 
For the orthogonal system  of functions $H_n^{(1)}(x)=\exp(-x^2/2)H_n(x)$, denote their
Hilbert transforms by
\beq
G_n(x)=\fr{1}{\pi}\int_{-\infty}^\infty \fr{H_n^{(1)}(t)dt}{t-x}, \quad n=0,1,\ldots, \quad -\infty<x<\infty.
\label{2.3}
\eeq
 Then the  functions $G_n(x)$ form a complete orthogonal system in the space $L_2(-\infty,\infty)$,
 \beq
\int_{-\infty}^\infty G_n(x)G_m(x)dx=\sqrt{\pi}2^nn!\Gd_{mn},
\label{2.4}
\eeq
and admit the following representations:
$$
G_{2m}(x)=-\sqrt{2}(2m)!x e^{-x^2/2}\sum_{k=0}^m\fr{(-1)^k 2^k\GF(1/2-k,3/2;x^2/2)}{(m-k)!\GG(k+1/2)},
$$
\beq
G_{2m+1}(x)=\sqrt{2}(2m+1)! e^{-x^2/2}\sum_{k=0}^m\fr{(-1)^k 2^k\GF(-1/2-k,1/2;x^2/2)}{(m-k)!\GG(k+3/2)}, \quad m=0,1,\ldots.
\label{2.11}
\eeq
Here,  $\Gd_{mn}$ is the Kronecker symbol,  $\GF(a,c;x)$ is the confluent hypergeometric function,
\beq
\GF(a,c;x)=\sum_{k=0}^\infty \fr{(a)_k x^k}{(c)_k k!},
\label{2.12}
\eeq
and $(a)_k$  is the factorial symbol, $(a)_k=a(a+1)\ldots(a+k-1)$. 

\vspace{2mm}

{\it Proof.} 
The first statement of the theorem follows from the unitarity of the Hilbert operator, the Parseval's relation
\beq
\int_{-\infty}^\infty G_n(x)G_m(x)dx=\int_{-\infty}^\infty H_n^{(1)}(x)H_m^{(1)}(x)dx,
\label{2.4'}
\eeq
and the orthogonality of the functions $H_n^{(1)}(x)$.
To prove formulas (\ref{2.11}), we 
 apply the Fourier transform to equation (\ref{2.3}) and  employ the convolution theorem
and the spectral relation for the Fourier operator (\cite{gra}, 7.376.1, p.804)
\beq
F[H^{(1)}_n(x)](\Gl)=i^n\sqrt{2\pi} H^{(1)}_n(\Gl), \quad n=0,1,\ldots,
\label{2.5}
\eeq
to deduce
\beq
F[G_n(x)](\Gl)=-i^{n+1}{\sgn}\Gl\sqrt{2\pi} H^{(1)}_n(\Gl), \quad n=0,1,\ldots.
\label{2.6}
\eeq
Here,
\beq
F[\Gf(x)](\Gl)=\int_{-\infty}^\infty\Gf(x)e^{i\Gl x}dx.
\label{2.7}
\eeq
The Fourier inversion applied for even and odd indices yields
the following alternative integral representations for the $G$-functions:
$$
G_{2m}(x)=(-1)^{m+1}\sqrt{\fr{2}{\pi}}
\int_0^\infty H^{(1)}_{2m}(\Gl)\sin \Gl x d\Gl,
$$
\beq
G_{2m+1}(x)=(-1)^{m}\sqrt{\fr{2}{\pi}}
\int_0^\infty H^{(1)}_{2m+1}(\Gl)\cos \Gl x d\Gl,\quad m=0,1,\ldots, \quad -\infty<x<\infty.
\label{2.8}
\eeq	
Express next the Hermite polynomials through the Laguerre polynomials  (\cite{bat2}, 10.13 (2), (3), p.193)
\beq
H_{2m}(x)=(-1)^m 2^{2m}m! L_m^{-1/2}(x^2), \quad H_{2m+1}(x)=(-1)^m 2^{2m+1}m! xL_m^{1/2}(x^2),
\label{2.9}
\eeq
where the Laguerre polynomials are given by (\cite{gra}, 8.970.1, p. 1000)
\beq
L_m^{\Ga}(x)=\sum_{k=0}^m \left(
\begin{array}{c}m+\Ga\\ m-k\\
\end{array}
\right)\fr{(-1)^k x^k}{k!}, \quad 
 \left(
\begin{array}{c} a\\ n\\
\end{array}
\right)=\fr{\GG(a+1)}{n! \GG(a-n+1)}.
\label{2.10}
\eeq
Substitute now the expressions (\ref{2.9}) into (\ref{2.8}) and use formula (\ref{2.10}) and the sine- and cosine-integral
transforms of the function $x^\Gb e^{-\Ga x^2}$ (\cite{bat3}, 2.4(24), p.74). After simple rearrangement we have ultimately the 
representations  (\ref{2.11}).

\vspace{2mm}

{\sc  Corollary 2.2} The semi-infinite Hilbert transforms of the weighted Laguerre polynomials $\Gn^{\mp1/2}e^{-\Gn/2}L_m^{-1/2}(\Gn)$
are given by
$$
\fr{1}{\pi}\int_0^\infty\fr{e^{-\Gn/2}L_m^{-1/2}(\Gn)d\Gn}{(\Gn-\Gx)\sqrt{\Gn}}=\fr{(-1)^mG_{2m}(\sqrt{\xi})}{2^{2m}m!\sqrt{\Gx}},
$$
\beq
\fr{1}{\pi}\int_0^\infty\fr{e^{-\Gn/2}L_m^{1/2}(\Gn)\sqrt{\Gn}d\Gn}{\Gn-\Gx}=\fr{(-1)^mG_{2m+1}(\sqrt{\xi})}{2^{2m+1}m!},\quad m=0,1,\ldots, \quad 0<\Gx<\infty.
\label{2.13}
\eeq

\vspace{2mm}

{\it Proof.}
Write down the Hilbert transforms (\ref{2.3}) separately for even and odd indices, make the substitutions
$\Gx=x^2$ and $\Gn=t^2$,  and employ formulas (\ref{2.9}). This brings us to the integral relations (\ref{2.13}).

\vspace{2mm}

The orthogonality relations (\ref{2.4}), when written for the functions $G_n(\sqrt{\Gx})$, imply 
$$
\int_0^\infty G_{2m}(\sqrt{\Gx})G_{2n}(\sqrt{\Gx})\fr{d\Gx}{\sqrt{\Gx}}=2^{2n}(2n)!\sqrt{\pi}\Gd_{mn},
$$
\beq
\int_0^\infty G_{2m+1}(\sqrt{\Gx})G_{2n+1}(\sqrt{\Gx})\fr{d\Gx}{\sqrt{\Gx}}=2^{2n+1}(2n+1)!\sqrt{\pi}\Gd_{mn}, \quad m,n=0,1,\ldots.
\label{2.14}
\eeq
Notice that the orthogonality relations (\ref{2.4}) can alternatively be derived by employing the orthogonality relation for the weighted Hermite polynomials $H_n^{(1)}(x)$
and formulas (\ref{2.8}). Indeed, for even indices we have
\beq
\int_{-\infty}^\infty G_{2n}(x)G_{2m}(x)dx=2(-1)^{m+1}\sqrt{\fr{2}{\pi}}\int_0^\infty G_{2n}(x)dx\int_0^\infty e^{-\Gl^2/2}H_{2m}(\Gl)\sin \Gl xd\Gl.
\label{2.15}
\eeq
On changing the order of integration and applying the inverse sine-transform
\beq
\int_0^\infty G_{2n}(x)\sin\Gl xdx=(-1)^{n+1}\sqrt{\fr{\pi}{2}}e^{-\Gl^2/2}H_{2n}(\Gl), \quad \Gl>0,
\label{2.16}
\eeq
we deduce formula (\ref{2.4}) for even indices. In the same fashion, this formula is derived for odd indices.

\subsection{Relations for the weighted Laguerre polynomials $e^{-\Gn}\Gn^{\pm 1/2}L_m(\Gn)$}

Consider the Hilbert transform of the Hermite polynomials with the new weight  $e^{-x^2}$, $H_n^{(2)}(x)=e^{-x^2}H_n(x)$.

\vspace{2mm}

{\sc Theorem 2.3} 
Denote the semi-infinite  Hilbert transform of the weighted Hermite polynomials $H_n^{(2)}(x)$ by
\beq
V_n(x)=\fr{1}{\pi}\int_{-\infty}^\infty \fr{H_n^{(2)}(t)dt}{t-x},  \quad n=0,1,\ldots, \quad -\infty<x<\infty.
\label{2.17}
\eeq
Then the functions $V_n(x)$ admit the following representations 
 in terms of the confluent hypergeometric function: 
$$
V_{2m}(x)=\fr{(-1)^{m+1}}{\sqrt{\pi}}2^{2m+1}m! xe^{-x^2}\GF\left(\fr12-m,\fr32; x^2\right),
$$
\beq
V_{2m+1}(x)=\fr{(-1)^{m}}{\sqrt{\pi}}2^{2m+1}m! e^{-x^2}\GF\left(-\fr12-m,\fr12; x^2\right),\;  m=0,1,\ldots, \; -\infty<x<\infty,
\label{2.21}
\eeq
 
 \vspace{2mm}

{\it Proof.} As before, apply the Fourier transform and use the convolution theorem to obtain
\beq
F[V_n(x)](\Gl)=-i\sgn\Gl F[H^{(2)}_n(x)](\Gl).
\label{2.18}
\eeq 
To compute the Fourier transform of the function $H^{(2)}_n(x)$, we consider the even and odd indices cases separately
and
employ the table integrals (\cite{gra}, 7.388.1, p. 806) 
$$
\int_0^\infty H^{(2)}_{2m}(t)\cos\Gl tdt=\fr12(-1)^m\sqrt{\pi}\Gl^{2m}e^{-\Gl^2/4},
$$
\beq
\int_0^\infty H^{(2)}_{2m+1}(t)\sin\Gl tdt=\fr12(-1)^m\sqrt{\pi}\Gl^{2m+1}e^{-\Gl^2/4}.
\label{2.19}
\eeq
If combined, these give
$$
F[V_{2m}(x)](\Gl)=-i \sqrt{\pi}(-1)^m\sgn\Gl \Gl^{2m} e^{-\Gl^2/4}, 
$$
\beq
F[V_{2m+1}(x)](\Gl)=\sqrt{\pi}(-1)^m|\Gl|^{2m+1} e^{-\Gl^2/4}, 
\quad m=0,1,\ldots, \quad -\infty<\Gl<\infty.
\label{2.20}
\eeq
By the Fourier inversion we deduce analogues of formulas (\ref{2.11})
and express the functions $V_n(x)$,  the Hilbert transforms of the functions $H_n^{(2)}(x)$,
 in terms of the confluent hypergeometric functions  by (\ref{2.21}).

\vspace{2mm}

On transforming 
the interval $(-\infty,\infty)$ into the semi-infinite interval and making the substitutions $x^2=\Gx$ and $t^2=\Gn$ 
formulas (\ref{2.17}) and (\ref{2.21})  enable us to prove the following result.

{\sc  Corollary 2.4} The semi-infinite Hilbert transforms of the weighted Laguerre polynomials $\Gn^{\mp1/2}e^{-\Gn}L_m^{-1/2}(\Gn)$
are expressed through the confluent hypergeometric functions  by 
$$
\fr{1}{\pi}\int_0^\infty\fr{e^{-\Gn}L_m^{-1/2}(\Gn)d\Gn}
{(\Gn-\Gx)\sqrt{\Gn}}=-\fr{2e^{-\Gx}}{\sqrt{\pi}}\GF\left(\fr12-m,\fr32; \Gx\right),
$$
\beq
\fr{1}{\pi}\int_0^\infty\fr{e^{-\Gn}L_m^{1/2}(\Gn)\sqrt{\Gn}d\Gn}{\Gn-\Gx}=
\fr{e^{-\Gx}}{\sqrt{\pi}}\GF\left(-\fr12-m,\fr12; \Gx\right),\quad m=0,1,\ldots, \quad 0<\Gx<\infty.
\label{2.22}
\eeq

\vspace{2mm}

Undoubtedly, the integral relations (\ref{2.22}) are simpler than (\ref{2.13}). However, if applied to singular integral equations with the Cauchy kernel in a semi-infinite interval,
they have a disadvantageous feature: the right-hand sides of the relations (\ref{2.22}), the functions
\beq
\Gf_m^{(1)}(\Gx)=-\fr{2e^{-\Gx}}{\sqrt{\pi}}
\GF\left(\fr12-m,\fr32; \Gx\right), \quad 
\Gf^{(2)}_m(\Gx)=\fr{e^{-\Gx}}{\sqrt{\pi}}\GF\left(-\fr12-m,\fr12; \Gx\right)
\label{2.23}
\eeq
do not form an orthogonal system in the space $L^2(0,\infty)$. 
At the same time, the systems $\{\Gf^{(1)}_m(\Gx)\}_{m=0}^\infty$ and $\{\Gf^{(2)}_m(\Gx)\}_{m=0}^\infty$ are linearly independent.
To prove the linear independence of the first system, denote 
\beq
\Gc_m(\Gn)=e^{-\Gn}\Gn^{-1/2}[c_0L_0^{-1/2}(\Gn)+c_1L_1^{-1/2}(\Gn)+\ldots+c_mL_m^{-1/2}(\Gn)].
\label{2.24}
\eeq
From the first formula in (\ref{2.22}) we deduce
\beq
\fr{1}{\pi}\int_0^\infty\fr{\Gc_m(\Gn)d\Gn}{\Gn-\Gx}=c_0\Gf_0^{(1)}(\Gx)+c_1\Gf_1^{(1)}(\Gx)+\ldots c_m\Gf_m^{(1)}(\Gx), \quad 0<\Gx<\infty.
\label{2.25}
\eeq
Suppose 
\beq
c_0\Gf_0^{(1)}(\Gx)+c_1\Gf_1^{(1)}(\Gx)+\ldots c_m\Gf_m^{(1)}(\Gx)=0, \quad 0<\Gx<\infty.
\label{2.26}
\eeq
Since the homogeneous singular integral equation 
\beq
\fr1{\pi}\int_0^\infty\fr{\Gc_m(\Gn)d\Gn}{\Gn-\Gx}=0,  \quad 0<\Gx<\infty,
\label{2.27}
\eeq
in the class of integrable in the interval $(0,\infty)$ functions
has the trivial solution only, we have
\beq
c_0L_0^{-1/2}(\Gn)+c_1L_1^{-1/2}(\Gn)+\ldots+c_mL_m^{-1/2}(\Gn)=0, \quad 0<\Gx<\infty.
\label{2.28}
\eeq
Now, the system of the Laguerre polynomials $\{L_n^\Ga(\Gn)\}_{n=0}^m$ is linearly independent in $(0,\infty)$. Therefore
$c_0=c_1=\ldots=c_m=0$, and the functions $\Gf_0^{(1)}(\Gx)$, $\Gf_0^{(1)}(\Gx)$, $\ldots,$ $\Gf_m^{(1)}(\Gx)$ are linearly independent for any $m$. 
We can show that the second system  $\{\Gf^{(2)}_m(\Gx)\}_{m=0}^\infty$  is linearly independent in a similar manner.

In what follows we orthogonalize the systems $\{\Gf^{(1)}_m(\Gx)\}_{m=0}^\infty$  and $\{\Gf^{(2)}_m(\Gx)\}_{m=0}^\infty$
and represent the elements of these orthogonal   systems as linear combinations of the functions $\Gf^{(1)}_m(\Gx)$ and  $\Gf^{(1)}_m(\Gx)$, respectively.
Consider the integrals
$$
J_m^{(1)}(\Gx)=\fr{1}{\pi}\int_0^\infty\fr{e^{-\Gn/2}L_m^{-1/2}(\Gn)d\Gn}{(\Gn-\Gx)\sqrt{\Gn}}, 
$$
\beq
J_m^{(2)}(\Gx)=\fr{1}{\pi}\int_0^\infty\fr{e^{-\Gn/2}L_m^{1/2}(\Gn)\sqrt{\Gn}d\Gn}{\Gn-\Gx},\quad m=0,1,\ldots, \quad 0<\Gx<\infty.
\label{2.29}
\eeq
We aim to compute them by utilizing formulas (\ref{2.22}). On making the substitutions $\Gn=2u$ and $\Gx=2t$ 
and employing the identity (\cite{bat2}, (40), p.192) that is
\beq
L_m^{-1/2}(2u)=\sum_{k=0}^m 
\left(
\begin{array}{c}
m-1/2\\ m-k\\
\end{array}
\right)(-1)^{m-k}2^kL_k^{-1/2}(u),
\label{2.30}
\eeq
we have for the integral $J_m^{(1)}(\Gx)$ 
\beq
J_m^{(1)}(2t)=\fr{1}{\sqrt{2}}\sum_{k=0}^m 
\left(
\begin{array}{c}
m-1/2\\ m-k\\
\end{array}
\right)\fr{(-1)^{m-k}2^k}{\pi}
\int_0^\infty\fr{e^{-u}L_k^{-1/2}(u)du}{(u-t)\sqrt{u}}.
\label{2.31}
\eeq 
The integral in (\ref{2.31}) is given by the first formula in  (\ref{2.22}).  This brings us to the following relation:
\beq
J_m^{(1)}(2t)=-\sqrt{\fr{2}{\pi}}\sum_{k=0}^m 
\left(
\begin{array}{c}
m-1/2\\ m-k\\
\end{array}
\right)(-1)^{m-k}2^ke^{-t}\GF\left(\fr12-k,\fr32; t\right).
\label{2.32}
\eeq
Analysis of this formula shows that 
the function $J_m^{(1)}(2t)$ is a linear combination of the functions $\Gf_k^{(1)}(t)$. Simple transformations yield ultimately
\beq
J_m^{(1)}(\Gx)=\fr{(-1)^mG_{2m}(\sqrt{\xi})}{2^{2m}m!\sqrt{\Gx}},\quad m=0,1,\ldots, \quad 0<\Gx<\infty.
\label{2.33}
\eeq
In a similar fashion we obtain 
\beq
J_m^{(2)}(\Gx)
=\fr{(-1)^mG_{2m+1}(\sqrt{\xi})}{2^{2m+1}m!},\quad m=0,1,\ldots, \quad 0<\Gx<\infty.
\label{2.34}
\eeq
Thus, we have deduced that the orthogonalization of the systems $\{\Gf^{(1)}_m(\Gx)\}_{m=0}^\infty$  and $\{\Gf^{(2)}_m(\Gx)\}_{m=0}^\infty$ 
leads to the integral relations (\ref{2.13}) derived in Section 2.1.

To complete this section, we invert the relations (\ref{2.22}) by representing them as the integral equation with the Cauchy kernel in a semi-infinite segment
\beq
\fr{1}{\pi}\int_0^\infty \fr{\Gf(\Gn)d\Gn}{\Gn-\Gx}=f(\Gx), \quad 0<\Gx<\infty.
\label{2.34.1}
\eeq
Its solution in the class of integrable functions  in $(0,\infty)$ and unbounded at $\Gx=0$
is
\beq
\Gf(\Gx)=-\fr{1}{\pi\sqrt{\xi}} \int_0^\infty \fr{\sqrt{\Gn}f(\Gn)d\Gn}{\Gn-\Gx},
\label{2.34.2}
\eeq
and the solution bounded at the point $\Gx=0$ has the form
\beq
\Gf(\Gx)=-\fr{\sqrt{\xi}}{\pi} \int_0^\infty \fr{f(\Gn)d\Gn}{\sqrt{\Gn}(\Gn-\Gx)}.
\label{2.34.3}
\eeq
Note that one of the ways to obtain formulas  (\ref{2.34.2}) and (\ref{2.34.3}) is to employ the solution of the singular integral 
equation in the segment $(a,b)$ \cite{gak},  put $a=0$ and pass to the limit $b\to\infty$.
Upon employing these expressions for the inverse operators and the relations (\ref{2.22}) we obtain, respectively,
$$
\fr{1}{\pi}\int_0^\infty \fr{\sqrt{\Gn} e^{-\Gn}\GF(1/2-m,3/2;\Gn)d\Gn}{\Gn-\Gx}=\fr{\sqrt{\pi}}{2}e^{-\Gx} L_m^{-1/2}(\Gx),
$$
\beq
\fr{1}{\pi}\int_0^\infty \fr{ e^{-\Gn}\GF(-1/2-m,1/2;\Gn)d\Gn}{\sqrt{\Gn}(\Gn-\Gx)}=-\sqrt{\pi}e^{-\Gx} L_m^{1/2}(\Gx),\quad
m=0,1,\ldots, \quad 0<\Gx<\infty.
\label{2.34.4}
\eeq

\setcounter{equation}{0}

 \section{Properties of the $G$-functions}
 
 We have proved that the Hilbert transforms of  the functions $\exp(-x^2/2)H_n(x)$, the  functions $G_n(x)$,
 form an orthogonal system in the space $L_2(-\infty,\infty)$ with the $L^2$-norm $||G_n(x)||=\pi^{1/4}2^{n/2}\sqrt{n!}$,
 while the systems of the functions $\{G_{2n}(\sqrt{\Gx})\}$ and 
 $\{G_{2n+1}(\sqrt{\Gx})\}$ ($n=0,1,\ldots$) are two orthogonal bases for the weighted space $L^2_w(0,\infty)$ with weight $w(\Gx)=\Gx^{-1/2}$. 
It was also deduced that, up to certain constant factors, the functions $\Gx^{-1/2}G_{2m}(\Gx)$ and  $G_{2m+1}(\Gx)$ are 
the semi-infinite Hilbert transforms of the weighted Laguerre polynomials  $e^{-\Gn/2}\Gn^{-1/2} L_m^{-1/2}(\Gn)$ and
$e^{-\Gn/2}\Gn^{1/2} L_m^{1/2}(\Gn)$, respectively. 
 
 In this section we aim to show that the functions $G_{2m}(\Gx)$ and $G_{2m+1}(\Gx)$ satisfy certain ordinary differential equations and 
 also to study their asymptotics for small and large $\Gx.$ 
 It is known (\cite{bat2}, (13), p.193) that the function $H_{2m}^{(1)}(\Gl)=e^{-\Gl^2/2}H_{2m}(\Gl)$ satisfies the differential equation
 \beq
 \left(\fr{d^2}{d\Gl^2}+4m+1-\Gl^2\right)H_{2m}^{(1)}(\Gl)=0.
 \label{2.35}
 \eeq
 Now, the function $G_{2m}(x)$, up to a factor, is the sine-transform of the function $H_{2m}^{(1)}(\Gl)$,
\beq
 G_{2m}(x)=(-1)^{m+1}\sqrt{\fr{2}{\pi}}
\int_0^\infty H^{(1)}_{2m}(\Gl)\sin \Gl x d\Gl.
\label{2.35'}
\eeq
 By multiplying equation (\ref{2.35}) by $(-1)^{m+1}\sqrt{2/\pi}\sin\Gl x$, integrating in $(0,\infty)$, and then integrating by parts 
 we deduce
 \beq
  \left(\fr{d^2}{dx^2}+4m+1-x^2\right)G_{2m}(x)=B^{(1)}_mx,\quad 0<x<\infty, \quad G_{2m}(0)=0,
 \label{2.36}
 \eeq
 where
 \beq
 B^{(1)}_m=\sqrt{\fr{2}{\pi}}\fr{(2m)!}{m!}.
 \label{2.37}
 \eeq
 In our derivations, we used the  fact that $H_{2m}^{(1)}(0)=(-1)^m(2m)!/m!$. 
 
 Similar actions discover that the $G$-functions of odd indices satisfy the differential
 equation
\beq
  \left(\fr{d^2}{dx^2}+4m+3-x^2\right)G_{2m+1}(x)=B^{(2)}_mx,\quad 0<x<\infty, \quad \fr{d}{dx}G_{2m+1}(0)=0,
 \label{2.38}
 \eeq
 where
 \beq
 B^{(2)}_m=\sqrt{\fr{2}{\pi}}\fr{2(2m+1)!}{m!}.
 \label{2.39}
 \eeq
 
Analyze now the asymptotics of the functions $G_n(\Gx)$ as $\Gx\to 0$ and $\Gx\to\infty$. 
The representations (\ref{2.11}) of $G_{2m}(x)$ and  $G_{2m+1}(x)$ and the series (\ref{2.12}) imply
\beq
G_{2m}(x)\sim a_0x,\quad G_{2m+1}(x)\sim a_1, \quad x\to 0,
\label{2.40}
\eeq
where
\beq
a_0=-\sqrt{2}(2m)!\sum_{k=0}^m\fr{(-1)^k 2^k}{(m-k)!\GG(k+1/2)},
\;
a_1=\sqrt{2}(2m+1)! \sum_{k=0}^m\fr{(-1)^k 2^k}{(m-k)!\GG(k+3/2)}.
\label{2.41}
\eeq
To derive the asymptotics of the $G$-functions for large $x$, we employ the formula (\cite{bat1}, 6.13.1(3), p. 278)
\beq
\GF(a,c;x)=\fr{\GG(c)}{\GG(a)}e^x x^{a-c}[1+O(x^{-1})], \quad x\to\infty.
\label{2.42}
\eeq
From the representations (\ref{2.11}) we deduce
$$
G_{2m}(x)\sim -\fr{2^{2m+1/2}\GG(m+1/2)}{\pi x},   \quad x\to\infty,
$$
\beq
G_{2m+1}(x)\sim -\fr{2^{2m+5/2}\GG(m+3/2)}{\pi x^2}, \quad x\to\infty. 
\label{2.43}
\eeq
It is possible to obtain full asymptotic expansions for large  $x$ for these functions by expressing the function $\GF$  in (\ref{2.11}) through
the Tricomi function $\Psi$ and  writing the asymptotic expansion of the function $\Psi$ (see \cite{bat1} 6.13.1(1) and
6.7(7), respectively). Alternatively, we may  use the relation (\ref{2.35'}) and the asymptotic  formula  (\cite{lig}, (3), p. 56)
\beq 
\int_0^\infty f(\Gl)\sin \Gl x d\Gl \sim \fr{f(0)}{x}-\fr{f''(0)}{x^3}+\fr{f^{IV}(0)}{x^5}-\ldots, \quad x\to\infty,
\label{2.43.1}
\eeq
This asymptotic expansion holds for all functions $f(x)$ defined with all its derivatives for $x\ge 0$. We discover for the $G$-functions of even indices
\beq
G_{2m}(x)\sim (-1)^{m+1}\sqrt{\fr{2}{\pi}}\left[\fr{H_m^{(1)}(0)}{x}-\fr{d^2 H_m^{(1)}(0)}{x^3dx^2}+\fr{d^4 H_m^{(1)}(0)}{x^5dx^4}-\ldots
\right],\quad x\to\infty,
\label{2.43.2}
\eeq
where
\beq
\fr{d^{2n} H_n^{(1)}(0)}{dx^{2n}}=
(-1)^{n+m}(2n)!(2m)!\sum_{j=0}^n\fr{2^{2j}}{(2j)!(m-j)!(2n-2j)!!},
\label{2.43.3}
\eeq
where  $(2m)!!=2\cdot 4\cdot\ldots\cdot (2m)$. On computing the first several terms we have
\beq
G_{2m}(x)\sim -\sqrt{\fr{2}{\pi}}\fr{(2m)!}{m!}\left(\fr{c_1}{x}+\fr{c_3}{x^3}+\fr{c_5}{x^5}+\fr{c_7}{x^7}+\ldots\right),\quad x\to\infty,
\label{2.43.4}
\eeq
where
\beq
c_1=1,\; c_3=4m+1, \; c_5=16m^2+8m+3,\;
c_7=64m^3+48m^2+68m+15.
\label{2.43.5}
\eeq
A similar asymptotic expansion may be obtained for  the functions $G_{2m+1}(x)$.
\begin{figure}[t]
\centerline{
\scalebox{0.7}{\includegraphics{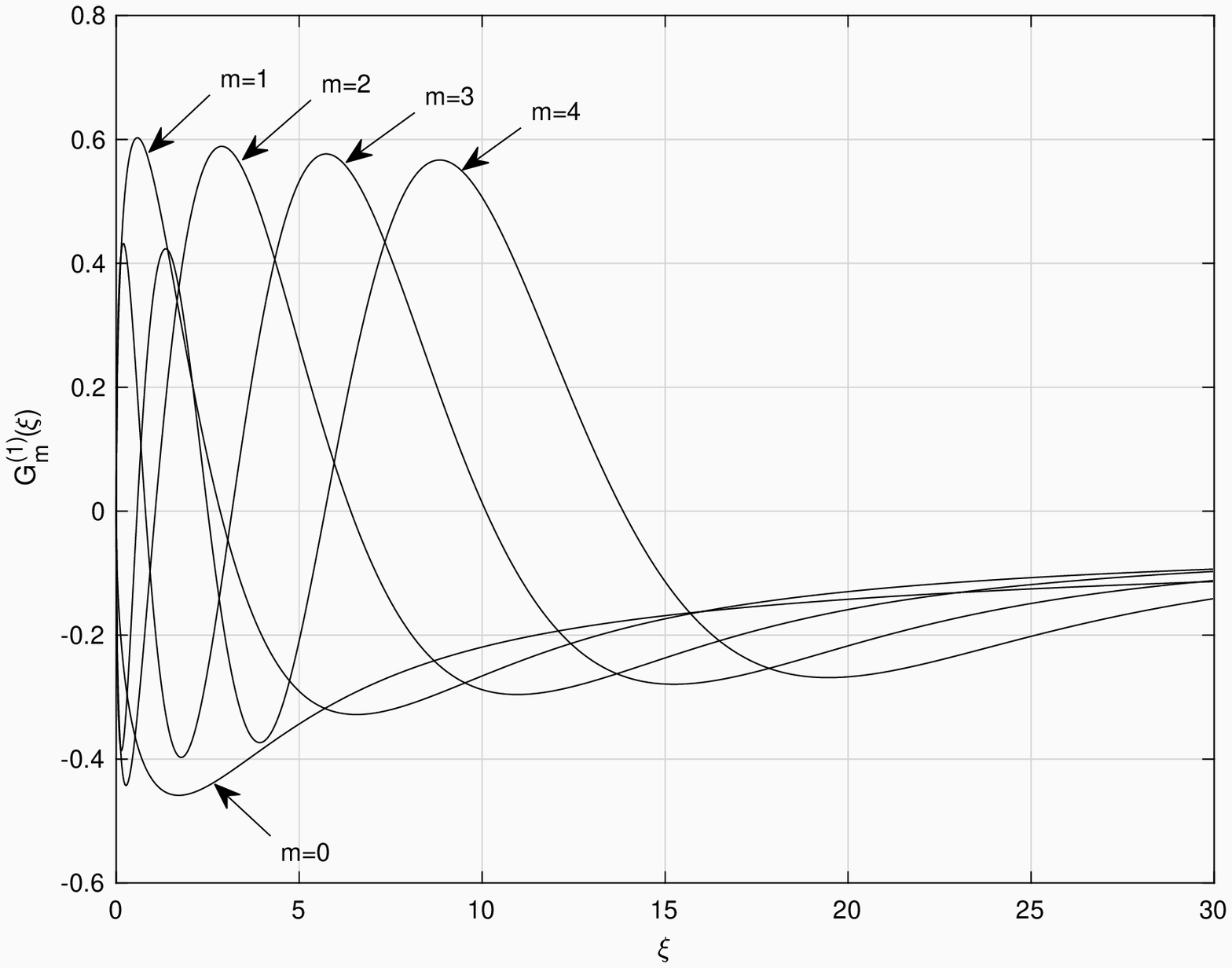}}}
\caption{The functions $G_m^{(1)}(\Gx)$, $m=0,1,2,3,4$.
}
\label{fig1}
\end{figure} 

\begin{figure}[t]
\centerline{
\scalebox{0.7}{\includegraphics{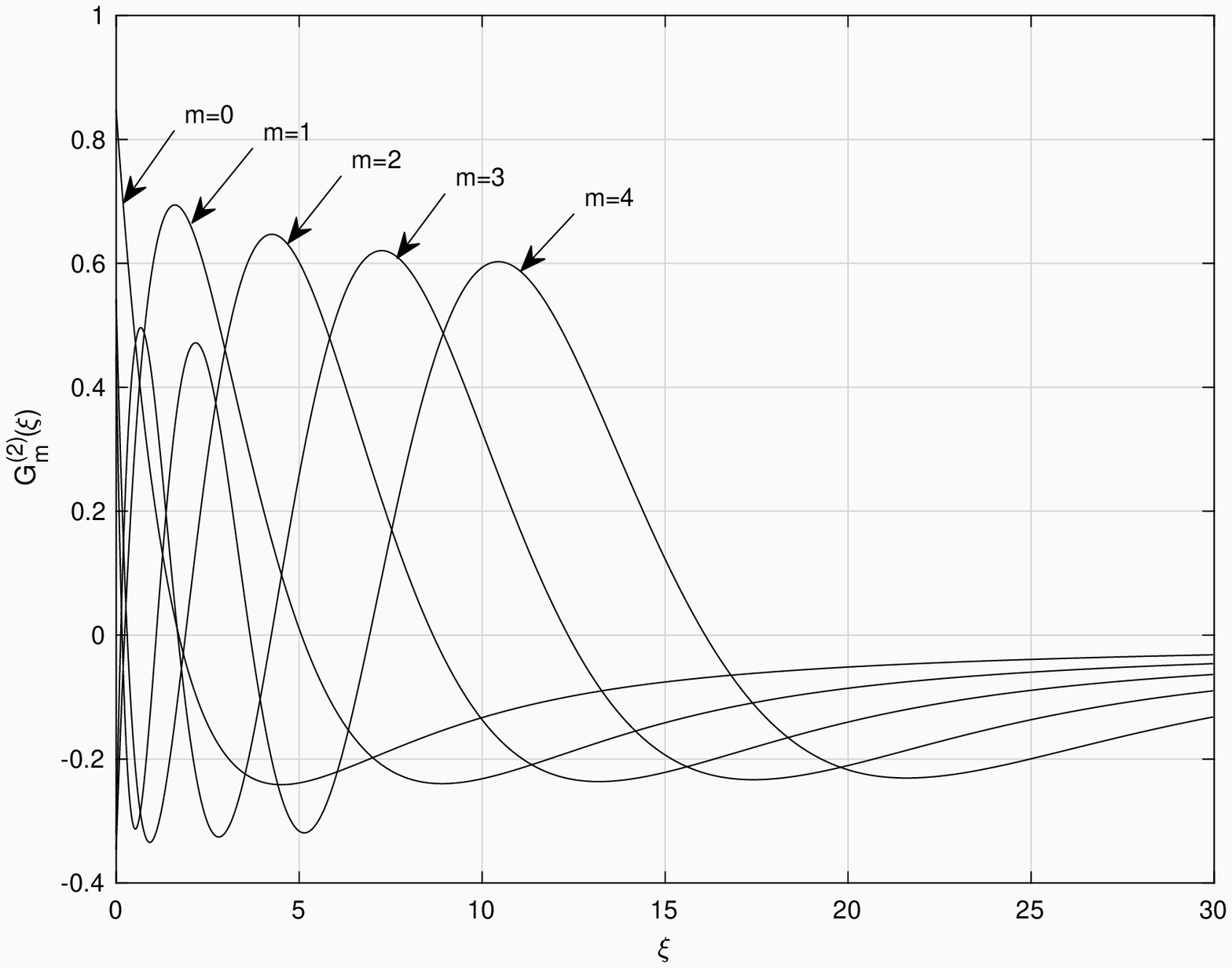}}}
\caption{The functions $G_m^{(2)}(\Gx)$, $m=0,1,2,3,4$.
}
\label{fig2}
\end{figure} 

For applications to integral equations, it will be convenient to denote 
$$
G_m^{(1)}(\xi)=\fr{G_{2m}(\sqrt{\Gx})}{2^m\sqrt{(2m)!}\pi^{1/4}},
$$
\beq
G_m^{(2)}(\xi)=\fr{G_{2m+1}(\sqrt{\Gx})}{2^{m+1/2}\sqrt{(2m+1)!}\pi^{1/4}},\quad m=0,1,\ldots, \quad 0<\xi<\infty.
\label{2.44}
\eeq
These functions form two orthonormal bases for the weighted space $L^2_w(0,\infty)$ with weight $w(\Gx)=\Gx^{-1/2}$,
\beq
\int_0^\infty G^{(j)}_{n}(\Gx)G^{(j)}_{m}(\Gx)\fr{d\Gx}{\sqrt{\Gx}}=\Gd_{mn}, \quad m,n=0,1,\ldots, \quad j=1,2.
\label{2.45}
\eeq
In Figures 1 and 2, we plot the
functions $G_m^{(1)}(\xi)$ and $G_m^{(2}(\xi)$ for $m=0,1,\ldots,4$, respectively. It has been discovered that although the functions
$G_m^{(1)}(\xi)$ and $G_m^{(2}(\xi)$ are not polynomials, in addition to the orthogonality and being a basis
of a certain weighted space, they share another property of classical orthogonal polynomials: in their interval of definition, $[0,\infty)$, the number of zeros correlates with
the index and equals 
 $m+1$ for both functions. Note that in the case of $G_m^{(1)}(\xi)$, $G_m^{(1)}(\xi)\sim \const \sqrt{\xi}$, $\Gx\to 0$, $m=0,1,\ldots$, and the point $\Gx=0$
 is counted as the first zero of the function $G_m^{(1)}(\xi)$. Due to the relations (\ref{2.43}) and $\ref{2.44})$, the functions $G_m^{(1)}(\xi)$ and $G_m^{(2}(\xi)$
 vanish at infinity,
 $$
 G_m^{(1)}(\Gx)\sim-\sqrt{\fr{2(2m-1)!!}{(2m)!!}}\fr{1}{\pi^{3/4}\sqrt{\Gx}},\quad
  G_m^{(2)}(\Gx)\sim-\sqrt{\fr{(2m+1)!!}{(2m)!!}}\fr{2}{ \pi^{3/4}\Gx},\quad \Gx\to\infty,
$$
where   $(2m\pm 1)!!=1\cdot 3\cdot\ldots\cdot (2m\pm 1)$.

\setcounter{equation}{0}
  
\section{Hilbert transforms associated with limiting relations for the Jacobi polynomials and
functions}\label{jac} 

\subsection{Laguerre polynomials and the Hilbert transform of the Jacobi polynomials}

A number of Hilbert transforms for special functions including the weighted Jacobi and Laguerre polynomials and  the confluent hypergeometric functions
can be derived from integral relations in a finite segment by letting a parameter involved go to infinity. To pursue this goal, we analyze the  
integral relation for the Jacobi polynomials \cite{tri1}
$$
\pi\cot\pi\Ga \,(1-x)^\Ga(1+x)^\Gb P_n^{(\Ga,\Gb)}(x)-\int_{-1}^1\fr{(1-t)^\Ga(1+t)^\Gb}{t-x}P_n^{(\Ga,\Gb)}(t)dt
$$
$$
=\fr{2^{\Ga+\Gb}\GG(\Ga)\GG(n+\Gb+1)}{\GG(n+\Ga+\Gb+1)}F\left(n+1,-n-\Ga-\Gb;1-\Ga; \fr{1-x}{2}\right), 
$$
\beq
 -1<x<1, \quad n=0,1,\ldots.
\label{4.1}
\eeq
Here, $\Ga>-1$, $\Ga\ne 0,1,\ldots$, $\Gb>-1$, $P_n^{(\Ga,\Gb)}(x)$ are the Jacobi polynomials, and $F(a,b;c;x)$ is the Gauss hypergeometric function.  
On making the substitutions $x=1-2\Gx/\Gb$ and $t=1-2\Gn/\Gb$ we infer 
$$
\pi\cot\pi\Ga \;\Gx^\Ga\left(1-\fr{\Gx}{\Gb}\right)^\Gb P_n^{(\Ga,\Gb)}\left(1-\fr{2\Gx}{\Gb}\right)+\int_{0}^\Gb\fr{\Gn^\Ga(1-\Gn/\Gb)^\Gb}{\Gn-\Gx}P_n^{(\Ga,\Gb)}\left(1-\fr{2\Gn}{\Gb}\right)d\Gn
$$
$$
=\fr{\Gb^\Ga \GG(\Ga)\GG(n+\Gb+1)}{\GG(n+\Ga+\Gb+1)}F\left(n+1,-n-\Ga-\Gb;1-\Ga; \fr{\Gx}{\Gb}\right), 
$$
\beq
 0<\Gx<\Gb, \quad n=0,1,\ldots.
\label{4.2}
\eeq
In what follows we use the connection between the Laguerre and  Jacobi polynomials (\cite{sze}, (5.3.4) p. 103)
\beq
L_n^\Ga(\Gx)=\lim_{\Gb\to\infty}P_n^{(\Ga,\Gb)}\left(1-\fr{2\Gx}{\Gb}\right)
\label{4.3}
\eeq
and the asymptotic formula for the $\GG$-functions
\beq
\fr{\GG(z+\Ga)}{\GG(z+\Gb)}\sim z^{\Ga-\Gb}, \quad z\to\infty.
\label{4.4}
\eeq
The last relation enables us to evaluate the limits
$$
\lim_{\Gb\to\infty}\fr{\Gb^\Ga\GG(n+\Gb+1)}{\GG(n+\Ga+\Gb+1)}=1,
$$
\beq
\lim_{\Gb\to\infty} F\left(n+1,-n-\Ga-\Gb;1-\Ga; \fr{\Gx}{\Gb}\right)=\GF(n+1,1-\Ga; -\Gx).
\label{4.5}
\eeq
Now, on passing to the limit $\Gb\to \infty$ in (\ref{4.2}) and using formula (\ref{4.3}) and Kummer's transformation (\cite{bat1}, (7), p. 253) 
\beq
\GF(n+1,1-\Ga; -\Gx)=e^{-\Gx}\GF(-\Ga-n,1-\Ga;\Gx),
\label{4.6}
\eeq
we obtain the following result.

\vspace{2mm}

{\sc Theorem 4.1} Let $\Ga>-1$, $\Ga\ne 0,1,\ldots$. Then
$$
\pi\cot\pi\Ga \;\Gx^\Ga e^{-\Gx} L_n^\Ga(\Gx)
+\int_{0}^\infty\fr{\Gn^\Ga e^{-\Gn} L_n^\Ga(\Gn)d\Gn}{\Gn-\Gx}
=\GG(\Ga)e^{-\Gx}\GF(-n-\Ga, 1-\Ga; \Gx), 
$$
\beq
0<\Gx<\infty, \quad n=0,1,\ldots.
\label{4.7}
\eeq

\vspace{2mm}

{\sc Remark 4.2} This theorem generalizes Corollary 2.4: formulas  (\ref{2.22}) can be immediately deduced from (\ref{4.7})  by putting there $\Ga=\pm 1/2$.

\subsection{Laguerre polynomials and  the Jacobi functions $Q_n^{(\Ga,\Gb)}(x)$ }

To derive an analogue of the integral relation (\ref{4.7}) for the interval $(-\infty,0)$, we analyze two representations
of the Jacobi function $Q_n^{(\Ga,\Gb)}(x)$ (\cite{sze}, (4.61.4), (4.61.5), p. 74)
\beq
Q_n^{(\Ga,\Gb)}(x)=-\fr{(x-1)^{-\Ga}(x+1)^{-\Gb}}{2}\int_{-1}^1(1-t)^\Ga(1+t)^\Gb\fr{P_n^{(\Ga,\Gb)}(t)dt}{t-x}, \quad n=0,1,\ldots,
\label{4.8}
\eeq
and
$$
Q_n^{(\Ga,\Gb)}(x)=\fr{2^{n+\Ga+\Gb}\GG(n+\Ga+1)\GG(n+\Gb+1)}{\GG(2n+\Ga+\Gb+2)}(x-1)^{-n-\Ga-1}(x+1)^{-\Gb}
$$
\beq
\times F\left(n+\Ga+1,n+1;2n+\Ga+\Gb
+2;\fr{2}{1-x}\right), \quad n=0,1,\ldots,
\label{4.9}
\eeq 
where $\Ga>-1$, $\Gb>-1$. The relations are valid in the whole complex plane cut along the segment $[-1,1]$.
We consider the case $x>1$, employ the variables $x=1-2\Gx/\Gb$ and $t=1-2\Gn/\Gb$, and denote
\beq
f_n(\Ga,\Gb;\Gx)=2(-\Gx)^\Ga\left(1-\fr{\Gx}{\Gb}\right)^\Gb Q_n^{(\Ga,\Gb)}\left(1-\fr{2\Gx}{\Gb}\right).
\label{4.10}
\eeq
This transforms the relations (\ref{4.8}) and (\ref{4.9}) to the following:
$$
f_n(\Ga,\Gb;\Gx)=\int_0^\Gb
\fr{\Gn^\Ga(1-\Gn/\Gb)^\Gb P_n^{(\Ga,\Gb)}(1-2\Gn/\Gb)d\Gn}{\Gn-\Gx}
$$
\beq
=\fr{\Gb^{n+\Ga+1}\GG(n+\Ga+1)\GG(n+\Gb+1)}{\GG(2n+\Ga+\Gb+2)(-\Gx)^{n+1}}
F\left(n+\Ga+1,n+1;2n+\Ga+\Gb+2;\fr{\Gb}{\Gx}\right).
\label{4.11}
\eeq 
Passing to the limit $\Gb\to\infty$ in the expression for the function $f_n(\Ga,\Gb;\Gx)$ gives 
\beq
\lim_{\Gb\to\infty} f_n(\Ga,\Gb;\Gx)=\GG(n+\Ga+1)(-\Gx)^{-n-1}
{}_2 F_0(n+1,n+\Ga+1;\Gx^{-1}), \quad \Gx<0,
\label{4.12}
\eeq
where ${}_2 F_0(\Ga,\Gb;z)$ is a generalized hypergeometric series that can also be expressed through the Tricomi function 
$\Psi$ (\cite{bat1}, 6.6(3), p. 257)
\beq
{}_2 F_0(\Ga,\Gb;-\Gx^{-1})=\Gx^\Ga\Psi(\Ga,\Ga-\Gb+1;\Gx).
\label{4.13}
\eeq
On letting $\Gb\to \infty$ in the first relation in (\ref{4.11}) in view of formula (\ref{4.3}), we have
 \beq
 \int_0^\infty\fr{\Gn^\Ga e^{-\Gn} L_n^\Ga(\Gn)d\Gn}{\Gn-\Gx}=\GG(n+\Ga+1)\Psi(n+1,1-\Ga;-\Gx), \quad n=0,1,\ldots; \quad -\infty<\Gx<0.
 \label{4.14}
 \eeq
This formula can be rewritten in terms of the the confluent hypergeometric function $\GF$ if we employ
the relation between the  $\Psi$- and $\GF$-functions (\cite{bat1}, 6.5(7), p.257). We have the following result.
\vspace{2mm}

{\sc Theorem 4.3} Let $\Ga>-1$, $\Ga\ne 0,1,\ldots$. Then
 $$
 \int_0^\infty\fr{\Gn^\Ga e^{-\Gn} L_n^\Ga(\Gn)d\Gn}{\Gn-\Gx}=\GG(\Ga)\GF(n+1,1-\Ga; -\Gx)
 $$
 \beq
+\fr{\GG(-\Ga)\GG(n+\Ga+1)}{n!}
 (-\Gx)^\Ga\GF(n+1+\Ga,1+\Ga;-\Gx), \quad n=0,1,\ldots, -\infty<\Gx<0.
 \label{4.15}
 \eeq

\vspace{2mm}

{\sc Remark 4.4} Alternatively, this relation can be derived by analytic continuation of the Gauss function in (\ref{4.11}) (\cite{bat1}, 2.10(2), p.108)
$$
F\left(n+\Ga+1,n+1;2n+\Ga+\Gb+2;\fr{\Gb}{\Gx}\right)=\fr{\GG(2n+\Ga+\Gb+2)\GG(\Ga)}{\GG(n+\Ga+1)\GG(n+\Ga+\Gb+1)}
$$
$$
\times \left(-\fr{\Gb}{\Gx}\right)^{-n-1}F\left(n+1,-n-\Ga-\Gb; 1-\Ga;\fr{\Gx}{\Gb}\right)+\fr{\GG(2n+\Ga+\Gb+2)\GG(-\Ga)}{\GG(n+1)\GG(n+\Gb+1)}
$$
\beq
\times\left(-\fr{\Gb}{\Gx}\right)^{-n-\Ga-1}F\left(n+\Ga+1,-n-\Gb; 1+\Ga;\fr{\Gx}{\Gb}\right), \quad \Gx<0,
\label{4.16}
\eeq
and consequent passing to the limit $\Gb\to\infty$ in (\ref{4.11}).

\subsection{Integral relations for cylindrical functions}

By passing to the limit $n\to\infty$ in certain relations for the Jacobi functions of the second kind it is possible to discover
some elegant formulas for cylindrical functions.

\vspace{2mm}

{\sc Theorem 4.5} Let $\Gl>0$ and $\Ga$ be a complex number such that $-1<\R\Ga<5/2$. Then
\beq
\fr{\pi}{\sin\pi\Ga}\left[z^{\Ga/2} J_{-\Ga}(\Gl\sqrt{z})-(-z)^{\Ga}z^{-\Ga/2}J_\Ga(\Gl\sqrt{z})\right]=\int_0^{\infty}
\fr{t^{\Ga/2}J_\Ga(\Gl\sqrt{t})dt}{t-z},
\label{4.23}
\eeq 
 where $J_\Ga(z)$ is the Bessel function, and $z^\Ga$ is the single branch in the $z$-plane cut along the ray $(-\infty,0]$ such that $\arg z\in[-\pi,\pi]$.

\vspace{2mm}

{\it Proof.}
We start with the following representation of the Jacobi function (\cite{bat2}, (19), p.171):
$$
Q_n^{(\Ga,\Gb)}(\Gz)=-\fr{\pi}{2\sin\pi\Ga}
P_n^{(\Ga,\Gb)}(\Gz)+\fr{2^{\Ga+\Gb-1}\GG(\Ga)\GG(n+\Gb+1)}{\GG(n+\Ga+\Gb+1)}(\Gz-1)^{-\Ga}(\Gz+1)^{-\Gb}
$$
\beq
\times F\left(n+1,-n-\Ga-\Gb;1-\Ga;\fr{1-\Gz}{2}\right) 
\label{4.17}
\eeq
valid in the whole $\Gz$-plane cut along the segment $[-1,1]$. On the cut sides, $\Gz=x\pm i0$, $\arg(\Gz-1)=\pm \pi$, $\arg(\Gz+1)=0$.
Put $\Gz=1-w^2/(2n^2)$. Then the branch cut $[-1,1]$ in the $\Gz$-plane is transformed into the cut $[0,2n]$ of the $w$-plane. 
Intending  to pass to the limit $n\to\infty$ we multiply equation (\ref{4.17}) by $n^{-\Ga}$ and use the limiting relation (\cite{bat2}, (41), p.173)
\beq
\lim_{n\to\infty} n^{-\Ga} P_n^{(\Ga,\Gb)}\left(1-\fr{w^2}{2n^2}\right)=\left(\fr{w}{2}\right)^{-\Ga} J_\Ga(w).
\label{4.18}
\eeq
This relation holds for arbitrary $\Ga$ and $\Gb$, uniformly in any bounded region of
the complex plane. It is directly verified that
$$
\lim_{n\to\infty}  F\left(n+1,-n-\Ga-\Gb;1-\Ga;\fr{w^2}{4n^2}\right) =\sum_{k=0}^\infty
\fr{(-1)^kw^{2k}}{k!(1-\Ga)_k 2^{2k}}
$$
\beq
=\GG(1-\Ga)\left(\fr{w}{2}\right)^{\Ga} J_{-\Ga}(w).
\label{4.19}
\eeq
Consequently,  we deduce from (\ref{4.17}) 
\beq
\lim_{n\to\infty}n^{-\Ga} Q_n^{(\Ga,\Gb)}\left(1-\fr{w^2}{2n^2}\right)=\fr{\pi 2^{\Ga-1}}{\sin\pi\Ga}\left[(-w^2)^{-\Ga} w^\Ga J_{-\Ga}(w)-w^{-\Ga}J_\Ga(w)
\right]. 
\label{4.20}
\eeq
Now, make the substitution $x=1-w^2/(2n^2)$ in the integral relation (\ref{4.8}) to obtain
$$
Q_n^{(\Ga,\Gb)}\left(1-\fr{w^2}{2n^2}\right)=(-w^2)^{-\Ga}\left(2-\fr{w^2}{2n^2}\right)^{-\Gb}
$$
\beq
\times\int_0^{2n}u^{2\Ga+1}\left(2-\fr{u^2}{2n^2}\right)^{\Gb}\fr{P_n^{(\Ga,\Gb)}(1-\fr12u^2/n^2)du}{u^2-w^2}.
\label{4.21}
\eeq
We multiply this equation by $n^{-\Ga}$, pass to the limit $n\to\infty$ and use formulas (\ref{4.18}) and  (\ref{4.20}). We have  
\beq
\fr{\pi}{2\sin\pi\Ga}\left[w^\Ga J_{-\Ga}(w)-(-w^2)^{\Ga}w^{-\Ga}J_\Ga(w)\right]=\int_0^{\infty}
\fr{u^{\Ga+1}J_\Ga(u)du}{u^2-w^2}.
\label{4.22}
\eeq
Next we make the substitutions $u=\Gl\sqrt{t}$ and $w=\Gl\sqrt{z}$, $\Gl$ is a positive parameter, and deduce ultimately formula (\ref{4.23}).

\vspace{2mm}

{\sc  Corollary 4.6} Let  $\Ga=-1/2-i\mu$ and $-\infty<\mu<\infty$. Then the Bessel function $J_\Ga(\Gl\sqrt{x}$ satisfies the integral relation
\beq
\fr{1}{\pi}\int_0^{\infty}
\fr{t^{\Ga/2}J_\Ga(\Gl\sqrt{t})dt}{t-x}+i\tanh\pi\mu x^{\Ga/2} J_\Ga(\Gl\sqrt{x})=
-\fr{x^{\Ga/2}}{\cosh\pi\mu}J_{-\Ga}(\Gl\sqrt{x}), \; 0<x<\infty.
\label{4.25}
\eeq

\vspace{2mm}

{\it Proof.}
Put $z=x\pm i0$, $0<x<\infty$, in the last relation. Since  $\arg(-z)=\mp \pi$, by the Sokhotski-Plemelj formulas
\beq
\fr{\pi x^{\Ga/2}}{\sin\pi\Ga}\left[J_{-\Ga}(\Gl\sqrt{x})-e^{\mp i\pi\Ga}J_\Ga(\Gl\sqrt{x})\right]=
\pm \pi i x^{\Ga/2}J_\Ga(\Gl\sqrt{x})+
\int_0^{\infty}
\fr{t^{\Ga/2}J_\Ga(\Gl\sqrt{t})dt}{t-x}.
\label{4.24}
\eeq
It is directly verified that both formulas may be put into the same form
(\ref{4.25}).

\vspace{2mm}

{\sc  Corollary 4.7}  Let $I_{\Ga}(x)$ be the modified Bessel function of the first kind,  $\Ga=-1/2-i\mu$,  and $-\infty<\mu<\infty$. Then
\beq
\fr{1}{\pi}\int_0^{\infty}
\fr{t^{\Ga/2}J_\Ga(\Gl\sqrt{t})dt}{t-x}=
\fr{(-x)^{\Ga/2}}{\cosh\pi\mu}[I_{\Ga}(\Gl\sqrt{-x})-I_{-\Ga}(\Gl\sqrt{-x})], \quad -\infty<x<0,
\label{4.30}
\eeq

\vspace{2mm}

{\it Proof.} Let $z\to x\pm i0$, $x<0$. Since $\arg z=\pi$, $\arg(-z)=0$, and
$I_{\Ga}(x)=e^{-i\pi\Ga/2}J_\Ga(ix)$, we deduce from (\ref{4.23}) the relation needed.

\vspace{2mm}

{\sc Remark 4.8}
On letting $\Gl\to 0^+$, we obtain from (\ref{4.24}) the following spectral relation for the operator $H$ in a semi-infinite interval:
\beq
\fr{1}{\pi i}\int_0^{\infty}
\fr{t^{-1/2+i\mu}dt}{t-x}=\tanh\pi\mu \,x^{-1/2+i\mu}, \quad 0<x<\infty, \quad -\infty<\mu<\infty.
\label{4.25'}
\eeq
This means that the function $f(t)=t^{-1/2+i\mu}$ is a generalized eigenfunction of the Hilbert operator in the interval $(0,\infty)$,
and $\tanh\pi\mu$ is its eigenvalue. We call $f(t)$ a generalized eigenfunction since it is not an $L_2(0,\infty)$-function.  In virtue of the inequality
$-\infty<\mu<\infty$,  the interval $(-1,1)$ is continuous spectrum of the operator $H$.

\vspace{2mm}

{\sc Remark 4.9} The particular case (\ref{4.25}) of the general formula (\ref{4.23}) can also be derived 
from the integral relation for the Jacobi polynomials (\cite{tri1}, \cite{pop})
$$
\fr{1}{\pi i}\int_{-1}^1\fr{P_n^{(\Ga,\bar\Ga)}(t)(1-t)^\Ga(1+t)^{\bar\Ga} dt}{t-x}-i\tanh\pi\mu P_n^{(\Ga,\bar\Ga)}(x)(1-x)^\Ga(1+x)^{\bar\Ga}
$$
\beq
=\left\{\begin{array}{cc}
(2\cosh\pi\mu)^{-1}P_{n-1}^{(\bar\Ga+1,\Ga+1)}(x), & n=1,2,\ldots,\\
0, & n=0,\\
\end{array}
\right., \quad -1<x<1,
\label{4.26}
\eeq
by utilizing the substitutions $x=1-\Gx^2/(2n^2)$ and $t=1-\Gn^2/(2n^2)$ and passing to the limit $n\to\infty$.

\vspace{2mm}

Finally, we show that the classical Hilbert relation \cite{tit}
\beq
\fr{1}{\pi}\int_{-\infty}^{\infty}\fr{\cos\Gl\Gn d\Gn}{\Gx-\Gn}=\sin\Gl\Gx, \quad -\infty<\Gx<\infty, \quad \Gl>0,
\label{4.29}
\eeq
can be deduced from (\ref{4.25}) as a particular case. Put 
$\mu=0$ in (\ref{4.25}). Due to the relations
\beq
J_{1/2}(z)=\sqrt{\fr{2}{\pi z}}\sin z, \quad J_{-1/2}(z)=\sqrt{\fr{2}{\pi z}}\cos z
\label{4.27}
\eeq
we immediately get
\beq
\fr{1}{\pi}\int_0^{\infty}\fr{\cos\Gl\sqrt{t} dt}{(x-t)\sqrt{t}}=\fr{\sin\Gl\sqrt{x}}{\sqrt{x}}, \quad 0<x<\infty.
\label{4.28}
\eeq
Equivalently, if the substitutions $\Gx=\sqrt{x}$ and $\Gn=\sqrt{t}$ are made, this may be written as the Hilbert relation (\ref{4.29}).

\setcounter{equation}{0}
\section{Applications to singular integral equations}\label{SIE} 

\subsection{Integral equation with the Cauchy kernel in a semi-infinite axis}

Based on the integral relations (\ref{2.13}) we derive an exact solution of the singular integral equation 
\beq
\fr{1}{\pi}\int_0^\infty \fr{\Gc(t)dt}{t-x}=f(x), \quad 0<x<\infty,
\label{5.1}
\eeq
in a series form free of singular integrals. In the class of functions  unbounded
at the point $x=0$, we expand the solution  through the Laguerre polynomials 
\beq
\Gc(x)=\fr{e^{-x/2}}{\sqrt{x}}\sum_{n=0}^\infty b_n \  L_n^{-1/2}(x).
\label{5.2}
\eeq
By substituting the series (\ref{5.2}) into equation (\ref{5.1}) and using the first formula in (\ref{2.13})  and the first orthogonality 
relation in (\ref{2.14}) we obtain for the coefficients $b_n$
\beq
b_n=\fr{(-1)^n n!}{\sqrt{\pi}(2n)!}\int_0^\infty f(t)G_{2n}(\sqrt{t})dt.
\label{5.3}
\eeq
In terms of the elements of the orthonormal basis $\{G_{n}^{(1)}(t)\}_{n=0}^\infty$ given by (\ref{2.44}) these coefficients have the form
\beq
b_n=\Ga_n\int_0^\infty f(t)G_{n}^{(1)}(t)dt,
\label{5.4}
\eeq
where 
\beq
 \Ga_n=\fr{(-1)^n}{\pi^{1/4}}\sqrt{\fr{(2n)!!}{(2n-1)!!}}.
\label{5.5}
\eeq

In the class of functions bounded at the point $x=0$ we seek the solution of equation (\ref{5.1}) in the form
\beq
\Gc(x)=e^{-x/2}\sqrt{x}\sum_{n=0}^\infty b_n \  L_n^{1/2}(x).
\label{5.6}
\eeq
On employing the second formulas in (\ref{2.13})  and  (\ref{2.14}) we derive the coefficient $b_n$ by quadratures possesing the $G$-functions of odd indices
\beq
b_n=\fr{(-1)^n n!}{\sqrt{\pi}(2n+1)!}\int_0^\infty f(t)G_{2n+1}(\sqrt{t})\fr{dt}{\sqrt{t}}.
\label{5.7}
\eeq
or, in terms of the orthonormal basis functions $G_n^{(2)}(x)$,
\beq
b_n=\Ga_n\sqrt{\fr{2}{2n+1}}
\int_0^\infty f(t)G_{n}^{(2)}(t)\fr{dt}{\sqrt{t}}.
\label{5.8}
\eeq

Interesting representations of the Cauchy kernel are derived by comparing the series- and closed-form solutions of the singular integral equation 
 (\ref{5.1}).  
 
\vspace{2mm}

{\sc Theorem 5.1} Let  $ 0<x<\infty$ and $ 0<t<\infty$. Then the following two bilinear expansions of the Cauchy kernel in terms of the Laguerre polynomials 
and the $G$-functions are valid:
$$
\fr{1}{t-x}=-\sqrt{\fr{\pi}{t}}e^{-x/2}\sum_{n=0}^\infty\fr{(-1)^n n!}{(2n)!}L_n^{-1/2}(x)G_{2n}(\sqrt{t}),
$$
\beq
\fr{1}{t-x}=-\sqrt{\pi}e^{-x/2}\sum_{n=0}^\infty\fr{(-1)^n n!}{(2n+1)!}L_n^{1/2}(x)G_{2n+1}(\sqrt{t}).
\label{5.8.0}
\eeq

\vspace{2mm}

{\it Proof.} The first representation of the Cauchy kernel is derived by comparing the series-form solution (\ref{5.2})  and its integral form (\ref{2.34.2}). 
Had we substituted the coefficients $b_n$ given by (\ref{5.7}) into the series (\ref{5.6}) and compared the new series with the 
integral-form solution (\ref{2.34.3}), we would have obtained the second formula.

\begin{figure}[t]
\centerline{
\scalebox{0.7}{\includegraphics{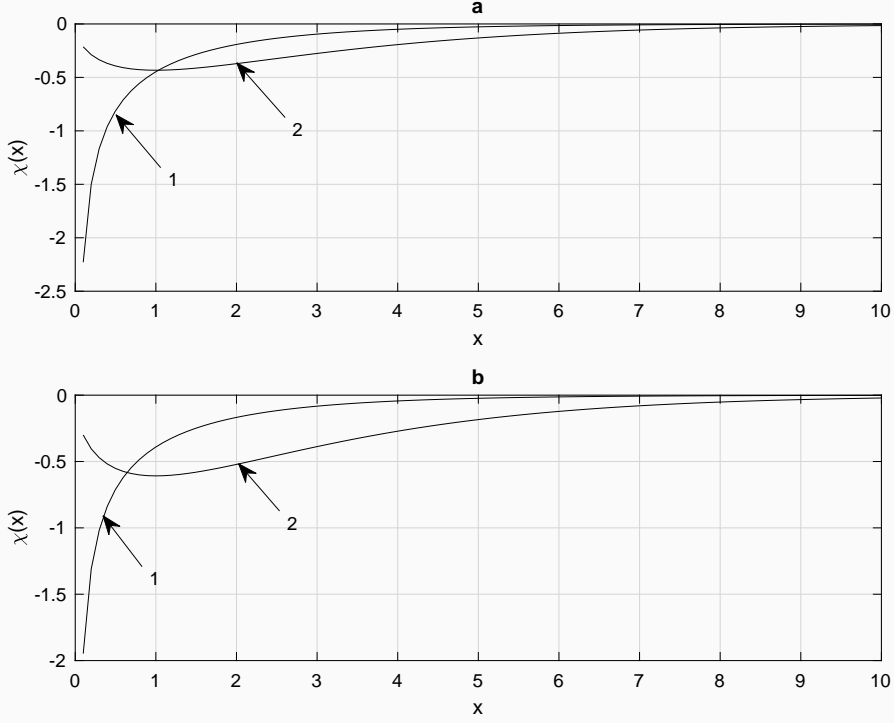}}}
\caption{The solution of the integral equation (\ref{5.1}) in the class of unbounded functions at $x=0$ (curve 1)
and bounded functions (curve 2). (a): $f(x)=\sqrt{x}e^{-x/2}$. (b): $f(x)=1$, $0<x<1$, and $f(x)=0$, $x>1$.} 
\label{fig3}
\end{figure}

\subsection{System of two singular integral equations}

Consider the systems of complete singular integral equations of the first kind with the Cauchy kernel
$$
\fr{1}{\pi}\int_0^\infty\left[\fr{1}{t-x}+K_{11}(t,x)\right]\Gc_1(t)dt+
\fr{1}{\pi}\int_0^\infty K_{12}(t,x) \Gc_2(t)dt=f_1(x), \quad 0<x<\infty,
$$
\beq
\fr{1}{\pi}\int_0^\infty K_{21}(t,x) \Gc_1(t)dt+\fr{1}{\pi}\int_0^\infty\left[\fr{1}{t-x}+K_{22}(t,x)\right]\Gc_2(t)dt=f_2(x), \quad 0<x<\infty.
\label{5.8.1}
\eeq
where $K_{jl}(t,x)$ may have a weak singularity at the line $x=t$. Suppose that the functions $\Gc_1(x)$ and $\Gc_2(x)$ are integrable in the interval $(0,\infty)$ and not bounded at the point $x=0$. Then, necessarily, 
they have the square root singularity at this point. We
expand the unknown functions in terms of the Laguerre polynomials 
\beq
\Gc_j(x)=\fr{e^{-x/2}}{\sqrt{x}}\sum_{m=0}^\infty b_m^{(j)} \  L_m^{-1/2}(x), \quad j=1,2.
\label{5.8.2}
\eeq
Because of  the kernels $k_{jl}$, in general, the coefficients $b_m^{(1)}$ and  $b_m^{(2)}$ cannot be found in explicit form.
By applying the same argument as in the case of the characteristic equation (\ref{5.1}) we deduce an infinite
system of linear algebraic equations
$$
b_n^{(1)}+\sum_{m=1}^\infty (c_{nm}^{(1,1)}b_m^{(1)}+ c_{nm}^{(1,2)}b_m^{(2)})=f_n^{(1)},
$$
\beq
b_n^{(2)}+\sum_{m=1}^\infty (c_{nm}^{(2,1)}b_m^{(1)}+ c_{nm}^{(2,2)}b_m^{(2)})=f_n^{(2)},
\quad n=0,1,\ldots.
\label{5.8.3}
\eeq
Here,
$$
c_{nm}^{(j,l)}=\fr{\Ga_n}{\pi}\int_0^\infty\int_0^\infty K_{jl}(t,x)\fr{e^{-t/2}}{\sqrt{t}}L_m^{-1/2}(t)G_n^{(1)}(x)dtdx, \quad j,l=1,2,
$$
\beq
f_n^{(j)}=\Ga_n\int_0^\infty f_j(x)G_n^{(1)}(x)dx, \quad j=1,2,
\label{5.9}
\eeq
where $\Ga_n$ are given by (\ref{5.5}). Assume that the kernels of the system of integral equations are chosen such that the system (\ref{5.8.3})
is regular.
By solving the infinite system (\ref{5.8.3}) by the reduction method we can approximately obtain the coefficients $b_m^{(j)}$ and therefore 
an approximate solution to the system  of singular integral equations.

\subsection{Integral relation for the Bessel function: a contact problem for a semi-infinite stamp}

Suppose a semi-infinite rigid stamp of profile $y=g(x)$ is indented into an elastic half-plane $|x|<\infty$, $-\infty<y<0$  such that the adhesion contact conditions
hold everywhere in the contact zone, while the rest of the boundary of the half-plane is free of traction,
$$
u(x,0)=c_1, \quad v(x,0)=g(x)+c_2, \quad 0<x<\infty,
$$
\beq
\Gs_y(x,0)=\tau_{xy}(x,0)=0, \quad -\infty<x<0.
\label{5.15}
\eeq
Here, $u$ and $v$ are the $x$- and $y$-components of the displacement vector, $\Gs_y$ and $\tau_{xy}$ are the stress tensor components, 
and $c_1$ and $c_2$ are constants.
Denote $p(x)=-\Gs_y(x,0)$ and $\tau(x)=-\tau_{xy}(x,0)$. Then this model problem is equivalent \cite{mus2}, \cite{gal} to the system
of integral equations
$$
\fr{\Gk-1}{\Gk+1}p(x)+\fr{1}{\pi}\int_0^\infty\fr{\tau(t)dt}{t-x}=0, \quad 0<x<\infty,
$$
\beq
\fr{\Gk-1}{\Gk+1}\tau(x)-\fr{1}{\pi}\int_0^\infty\fr{p(t)dt}{t-x}=\fr{4Gg'(x)}{\Gk+1}, \quad 0<x<\infty,
\label{5.16}
\eeq
where $\Gk=3-4\nu$, $\nu$ is the Poisson ratio, and $G$ is the shear modulus.
In terms of the function $\Gvf(x)=p(x)+i\tau(x)$, this system may be written as a single integral equation
\beq
\fr{1}{\pi}\int_0^\infty\fr{\Gvf(t)dt}{t-x}+i\tanh\pi\mu \, \Gvf(x)=f(x),\quad 0<x<\infty,
\label{5.17}
\eeq
where  $f(x)=-4G(\Gk+1)^{-1}g'(x)$ and $\mu=(2\pi)^{-1}\ln (3-4\nu)$.  For materials with the Poisson ratio $\nu\in(0,1/2)$, $\mu\in(0, \mu_0)$, $\mu_0\approx 0.17484958$.
This equation can exactly be solved by the Mellin transform or by the method of the Riemann-Hillbert
problem. In what follows, we propose  an alternative technique based on the integral relation for the Bessel function (\ref{4.25}) and the Hankel transform.
Represent the unknown function $\Gvf(x)$ in the integral form
\beq
\Gvf(x)=x^{\Ga/2}\int_0^\infty \Gc(\Gl)
J_\Ga(\Gl\sqrt{x})d\Gl, \quad \Ga=-\fr12-i\mu.
\label{5.18}
\eeq
Here, the density $\Gc(\Gl)$ is to be determined.  On substituting this integral into equation (\ref{5.17}), changing the order of integration, and
employing the  relation (\ref{4.25}) we 
deduce
\beq
\int_0^\infty \Gc(\Gl)J_{-\Ga}(\Gl\sqrt{x})d\Gl=-x^{-\Ga/2}\cosh\,\pi\mu \,f(x), \quad 0<x<\infty.
\label{5.19}
\eeq
By applying Hankel inversion we find the function $\Gc(\Gl)$ 
\beq
\Gc(\Gl)=-\fr{\Gl}{2}\cosh\pi\mu\int_0^\infty f(x)x^{-\Ga/2} J_{-\Ga}(\Gl\sqrt{x})dx.
\label{5.20}
\eeq
We remark that the profile $g(x)$ of the stamp  is assumed to be chosen such that the function $f(x)=-4G(\Gk+1)^{-1}g'(x)$ decays at infinity at the rate
sufficient for  the integral in (\ref{5.20}) being convergent.
It is directly verified that at the point $x=0$ and at infinity, the solution (\ref{5.18})  has the asymptotics required: it is oscillates and 
$\Gvf(x)=O(x^{-1/2})$, $x\to 0$,
and $\Gvf(x)=O(x^{-1/2})$, $x\to\infty$.

\setcounter{equation}{0}
\section{Quadrature formula for the Cauchy integral in  a semi-infinite interval}

In this section we obtain a quadrature formula for the Cauchy principal value of the singular integral
\beq
I^\Ga[f](x)=\fr{1}{\pi}\int_0^\infty \fr{f(t)w(t)dt}{t-x}, \quad w(t)=t^\Ga e^{-t}, \; \Ga>-1, \; \Ga\ne 0,1,\ldots, \; 0<x<\infty.
\label{6.1}
\eeq

\vspace{2mm}

{\sc Theorem 6.1} Let $f(x)$ be H\"older-continuous in any finite interval $[0,a]$, $a>0$, $|f(x)|\le Ce^{\Gb x}$, $x\to\infty$, $C=\const$, $\Gb<1$,    $w(t)=t^\Ga e^{-t}$, 
$\Ga>-1$, and $\Ga\ne 0,1,\ldots.$ Then 
\beq
I^\Ga[f](x)=\sum_{m=1}^n\Gg_m f(x_m)\fr{Q^\Ga_n(x)-Q^\Ga_n(x_m)}{x-x_m}+R_{n}(x), \quad 0<x<\infty,\quad x\ne x_m,
\label{6.13}
\eeq
where  
\beq
\Gg_m=-\fr{x_m}{(n+\Ga)L_{n-1}^\Ga(x_m)},
\label{6.13'}
\eeq
$x_m$ ($m=1,2,\ldots,n$) are the zeros of the degree-$n$ Laguerre polynomial $L_n^\Ga(x)$, and
\beq
Q^\Ga_n(x)=\fr{\GG(\Ga)}{\pi}e^{-x}\GF(-n-\Ga,1-\Ga;x)-\cot\pi\Ga x^\Ga e^{-x}L_n^\Ga(x).
\label{6.14}
\eeq
For $x=x_j$, 
\beq
I^\Ga[f](x_j)=\Gg_j f(x_j)\fr{dQ^\Ga_n(x_j)}{dx}+\sum_{m=1, m\ne j}^n\Gg_m f(x_m)\fr{Q^\Ga_n(x_j)-Q^\Ga_n(x_m)}{x_j-x_m}+R_{n}(x_j),
\label{6.15}
\eeq
where
$$
\fr{dQ^\Ga_n(x_j)}{dx}=\fr{1}{\pi}e^{-x_j}\GG(\Ga)\left[\fr{n+\Ga}{\Ga-1}\GF(-n-\Ga+1,2-\Ga;x_j)-\GF(-n-\Ga,1-\Ga;x_j)\right]
$$
\beq
+\cot\pi\Ga x_j^{\Ga-1}e^{-x_j}(n+\Ga)L_{n-1}^\Ga(x_j).
\label{6.16}
\eeq
Formula (\ref{6.13}) is exact, and $R_{n}(x)\equiv 0$ when $f(x)$ is a polynomial of degree not higher than $n-1$. Otherwise, if
$f(x)=M_{n-1}(x)+r(x)$ and $M_{n-1}(x)$ is  a polynomial of degree $n-1$, then
the reminder $R_{n}(x)$ of the quadrature formula (\ref{6.13}) is given by 
\beq
R_{n}(x)=-\sum_{m=1}^n\fr{\Gg_mr(x_m)}{x-x_m}[Q_n^\Ga(x)-Q_n^\Ga(x_m)]+\fr{1}{\pi}\int_0^\infty \fr{r(t)w(t)dt}{t-x}.
\label{6.12.7}
\eeq

\vspace{2mm}

{\it Proof.}
For our derivations, we use the method \cite{kor} proposed for the Cauchy integral in the finite segment
$(-1,1)$ and the integral relation (\ref{4.7})  of Theorem 4.1.
Introduce the system of orthonormal Laguerre polynomials
\beq
p_j(t)=h_j^{-1/2}L_j^\Ga(t), \quad h_j=\fr{\GG(\Ga+j+1)}{j!},\quad j=0,1,\ldots, 
\label{6.2}
\eeq
and assume first that $f(x)=M_{n-1}(x)$ is a polynomial of degree $n-1$. It will be convenient to express it in terms of the polynomials $p_j(t)$
\beq
f(t)=\sum_{j=0}^{n-1} f_j p_j(t),
\label{6.3}
\eeq
where
\beq
f_j=\int_0^\infty f(t)p_j(t)w(t)dt, \quad j=0,1,\ldots, n-1.
\label{6.4}
\eeq
By using the Gauss quadrature formula exact for polynomials of degree not higher than $2n-1$ we find
\beq
f_j=\sum_{m=1}^{n-1} A_m f(x_m)p_j(x_m)+\hat R_n(f),\quad j=0,1,\ldots,n-1,
\label{6.5}
\eeq
where $x_m$ ($m=1,2,\ldots,n$) are the zeros of the Laguerre polynomial $L_n^\Ga(x)$, $A_m$ are the Christoffel coefficients,
\beq
A_m=\fr{\GG(\Ga+n+1)}{n! x_m[\fr{d}{dx}L_n^\Ga(x_m)]^2},
\label{6.6}
\eeq
and $\hat R(f)$ is the reminder. Formula (\ref{6.5}) is exact when $f(t)=M_{n-1}(t)$. Denote further
\beq
q_j(x)=\int_0^\infty \fr{p_j(t)w(t)dt}{t-x}, \quad 0<x<\infty.
\label{6.6'}
\eeq
Substitute  the sum (\ref{6.3}) into (\ref{6.1}) and, in view of (\ref{6.5}) and (\ref{6.6}),  obtain for the principal part of the integral (\ref{6.1})
\beq
I^\Ga[f](x)=\fr{1}{\pi}\sum_{m=1}^n A_m f(x_m)\sum_{j=0}^{n-1} p_j(x_m)q_j(x).
\label{6.7}
\eeq
Write down next the Christoffel-Daurboux formula (\cite{sze}, (3.2.3), p. 43)
\beq
\sum_{j=0}^{n-1} p_j(x)p_j(t)=\fr{k_{n-1}}{k_n}\fr{p_n(x)p_{n-1}(t)-p_{n-1}(x)p_n(t)}{x-t},
\label{6.8}
\eeq
where $k_n=(-1)^n(n!\sqrt{h_n})^{-1}$. Integration of this identity with the weight $w(t)$ over the interval $(0,\infty)$ yields
\beq
p_{n-1}(x)q_n(x)-p_n(x)q_{n-1}(x)=\fr{k_n}{k_{n-1}}.
\label{6.8'}
\eeq
In particular,
\beq
q_n(x_m)=\fr{k_n}{k_{n-1}p_{n-1}(x_m)}.
\label{6.8''}
\eeq
On combining formulas (\ref{6.6'}), (\ref{6.8}), and (\ref{6.8'})
it is possible to establish the following identity \cite{kor}:
\beq
\sum_{j=0}^{n-1} p_j(x)q_j(t)
=\fr{k_{n-1}}{k_n(t-x)}
[p_n(x)(q_{n-1}(x)-q_{n-1}(t))-p_{n-1}(x))(q_{n}(x)-q_{n}(t))].
\label{6.9}
\eeq
Therefore, the internal sum in formula (\ref{6.7}) transforms to
\beq
\sum_{j=0}^{n-1} p_j(x_m)q_j(x)
=\fr{k_{n-1}p_{n-1}(x_m)}{k_n(x-x_m)}
[q_{n}(x)-q_{n}(x_m)].
\label{6.10}
\eeq
Now, by substituting this expression into formula (\ref{6.7}) we deduce
\beq
I^\Ga[f](x)=\fr{k_{n-1}}{\pi k_n}\sum_{m=1}^nA_m f(x_m)\fr{p_{n-1}(x_m)[q_n(x)-q_n(x_m)]}{x-x_m}.
\label{6.11}
\eeq
Finally, since $Q^\Ga_n(x)=\pi^{-1}\sqrt{h_n}q_n(x)$ and
\beq
\fr{k_{n-1}A_m}{k_n\sqrt{h_n h_{n-1}}}=-\fr{x_m}{(n+\Ga)[L_{n-1}^\Ga(x_m)]^2},
\label{6.12}
\eeq
we derive the quadrature formula (\ref{6.13}) for the singular integral (\ref{6.1}). Here, we employed an alternative
formula for the Christoffel coefficients (\ref{6.6})
\beq
A_m=\fr{\GG(\Ga+n)x_m}{n!(n+\Ga)[L_{n-1}^\Ga(x_m)]^2}.
\label{6.12'}
\eeq 
For $x=x_m$, by the L'H\^opital's rule we transform formula (\ref{6.13}) into the form (\ref{6.15}) with $R_{n}(x)=0$.

In the case $f(x)=M_{n-1}(x)+r(x)$ we derive the representation (\ref{6.12.7}) from (\ref{6.1}) and (\ref{6.13}).

\vspace{2mm}

{\sc  Corollary 6.2}  Let $x=\Gx_j$ be a zero of the function $Q_n^\Ga(x)$ given by (\ref{6.14}). Then the quadrature formula
(\ref{6.13}) has the form
\beq
I^\Ga[f](\Gx_j)=\fr{1}{\pi}\sum_{m=1}^n \fr{A_m f(x_m)}{x_m-\Gx_j}+R_n(\Gx_j),
\label{6.12.1}
\eeq
where $A_m$ are the Christoffel coefficients given by (\ref{6.12'}). This formula is exact for any polynomial $f(x)$ of degree $2n$.

\vspace{2mm}

{\it Proof.} In consequence of the relation (\ref{6.8''}) we obtain
\beq
Q_n^\Ga(x_m)=-\fr{\GG(\Ga+n)}{\pi n! L_{n-1}^\Ga(x_m)}.
\label{6.12.2}
\eeq
Putting $x=\Gx_j$ in (\ref{6.13}) and in view of  $Q_n^\Ga(\Gx_j)=0$ we have formula (\ref{6.12.1}). Since it is the Gauss quadrature formula 
in the interval $(0,\infty)$ associated with the Laguerre polynomials $L_n^\Ga(x)$, it is exact  for any polynomial $f(x)/(x-\Gx_j)$ of degree $2n-1$. Therefore
the last statement of Corollary 6.2 follows.

\vspace{2mm}

\begin{figure}[t]
\centerline{
\scalebox{0.7}{\includegraphics{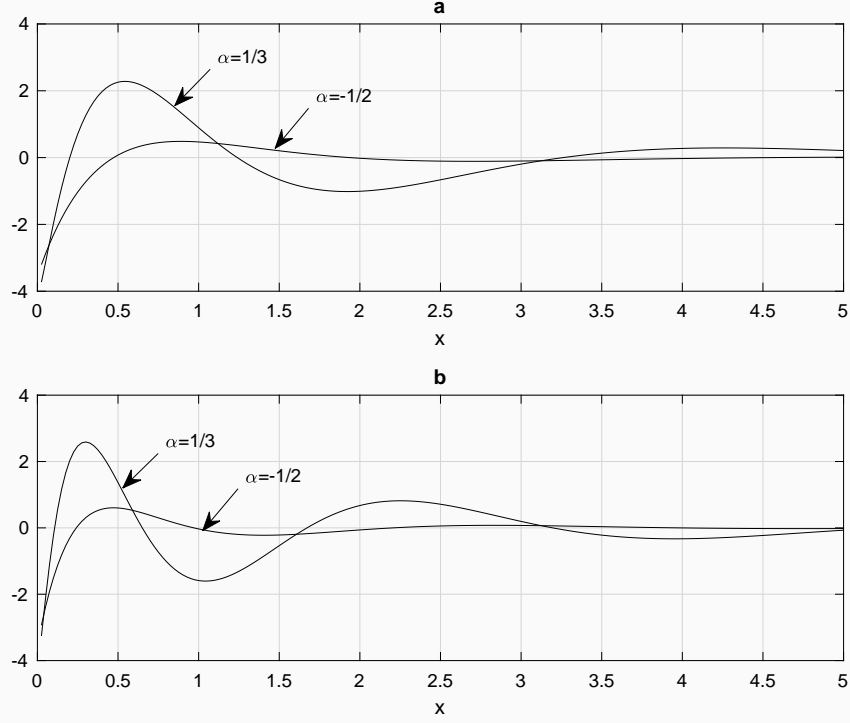}}}
\caption{The function $Q_n^\Ga(x)$ for $\Ga=1/3$ and $\Ga=-1/2$. (a): $n=5$. (b): $n=10$.
}
\label{fig4}
\end{figure}

\vspace{2mm}

{\sc Remark 6.3} Our numerical tests reveal that the function $Q_n^\Ga(x)$ has exactly $n+1$ zeros in the interval $(0,\infty)$ for any values of the parameters $n$ and 
$\Ga\in(-1,0)\cup(0,1)$ used.
Sample curves of the function $Q_n^\Ga(x)$ for $n=5$ and $n=10$ when $\Ga=1/3$ and $\Ga=-1/2$ are given in Figures 4a and 4b. It is also found that the amplitude of the function $Q_n^\Ga(x)$ rapidly decreases as $x\to\infty$.

\vspace{2mm}

Considering that a continuous function cannot be uniformly approximated in the interval $[0,\infty)$ by a polynomial, it is infeasible to find
an upper bound for the reminder in formula (\ref{6.15}) in the form $|R_{n}(x)|<\Gd_{n}$, $0\le x<\infty$, $\Gd_{n}\to 0$, $n\to\infty$,
for the general class of functions employed in Theorem 6.1. 
However, it is possible to estimate $R_{n}(x)$ for functions decaying at infinity as $x^{-m}$, $m>0$. 
\vspace{2mm}

{\sc Theorem 6.4}  Let  $f(x)$ be a continuously differentiable function in any finite segment  $[0,a]$, $a>0$, and
$f(x)\sim Cx^{-m}$, $x\to\infty$, $m>0$. For any $\Gve>0$ define a function
\beq
g_b(x)=\left\{\begin{array}{cc}
f'(x), & 0\le x\le b,\\
0, & x>b,\\
\end{array}
\right.
\quad b>\sqrt[m+1]{\fr{m|C|}{\Gve}}.
\label{8.12.2'}
\eeq
Let $P'_{n-1}(x)$ be the polynomial
of best approximation for the function $g_b(x)$ in the segment $[0,b]$ and
\beq
\max_{0\le x\le b}|g_b(x)-P'_{n-1}(x)|=e_{n}(b),
\label{6.12.4}
\eeq
Then the reminder $R_{n}(x)$ of the quadrature formula (\ref{6.13}) is estimated by
$$
|R_{n}(x)| \le 2\left(\fr{2\GG(\Ga+1)}{\pi}+|Q_n^\Ga(x)|\sum_{m=1}^n |\Gg_m|\right)\tilde e_{n}(b)
$$
\beq
\le 2\GG(\Ga+1)\left(\fr{2}{\pi}+\fr{n!|Q_n^\Ga(x)|}{\GG(n+\Ga)}\max_{1\le m\le n}|L_{n-1}^\Ga(x_m)|\right)\tilde e_{n}(b), \quad 0\le x<\infty,
\label{6.12.5}
\eeq
where $\tilde e_{n}(b)=\max\{\Gve, e_{n}(b)\}$. If $x=\Gx_j$, then the reminder has the bound
\beq
|R_n(\Gx_j)| \le\fr{4\GG(\Ga+1)}{\pi}\tilde e_{n}(b).
\label{6.12.6}
\eeq

\vspace{2mm}

{\it Proof.} Similar to \cite{kor}, rewrite formula (\ref{6.12.7}) for the reminder  in the form
\beq
R_n(x)=\sum_{m=1}^n\fr{\Gg_m[r(x)-r(x_m)]}{x-x_m}[Q_n^\Ga(x)-Q_n^\Ga(x_m)]+\fr{1}{\pi}\int_0^\infty \fr{r(t)-r(x)}{t-x}w(t)dt.
\label{6.12.8}
\eeq
By recalling that
\beq
\Gg_m Q_n^\Ga(x_m)=\fr{A_m}{\pi},
\label{6.12.9}
\eeq
we arrive at the formula
 \beq
R_n(x)=\sum_{m=1}^n\left[\Gg_mQ_n^\Ga(x)-\fr{A_m}{\pi}\right]
\fr{r(x)-r(x_m)}{x-x_m}+\fr{1}{\pi}\int_0^\infty\fr{r(t)-r(x)}{t-x}w(t)dt.
\label{6.12.3}
\eeq
Now, if we take into account formula (\ref{6.12.4}), the inequalities
\beq
|r'(x)|<|f'(x)-g_b(x)|+|g_b(x)-P'_{n-1}(x)|<\Gve+e_n(b),\quad 0\le x<\infty,
\label{6.12.9'}
\eeq
and also the integral
\beq
\int_0^\infty w(t)dt=\GG(\Ga+1),
\label{6.12.10}
\eeq
we deduce the first bound in (\ref{6.12.5}). To derive the second bound, we notice that $\Gg_mQ_n^\Ga(x_m)>0$,
\beq
\sum_{m=1}^n \Gg_mQ_n^\Ga(x_m)=\fr{1}{\pi}\int_0^\infty w(t)dt,
\label{6.12.11}
\eeq
and that, due to formula (\ref{6.8''}),
\beq
\fr{1}{Q_n^\Ga(x_m)}=\fr{\pi n!}{\GG(n+\Ga)}L_{n-1}^\Ga(x_m).
\label{6.12.12}
\eeq
Since $Q_n^\Ga(\Gx_j)=0$, the  bound (\ref{6.12.6}) follows from (\ref{6.12.5})  immediately.

\vspace{2mm}

Numerical tests implemented confirm numerical efficiency of the quadrature formula (\ref{6.13}). Sample curves for the integral $I^\Ga[f](x)$ for $\Ga=-1/2$ when the degree 
of the Laguerre polynomial $L_n^{-1/2}$is $n=5$ and $n=10$ are shown in Figures 5a -- 5d. 
In the case of Figure 5a, the integral can be evaluated exactly,
\beq
I^{-1/2}[\sqrt{t}](x)=\fr{1}{\pi}\int_0^\infty \fr{e^{-t}dt}{t-x}=-\fr{e^{-x}}{\pi}{\rm Ei}(x),\quad 0<x<\infty,
\label{6.17}
\eeq
where
\beq
{\rm Ei}(x)=\Gg+\ln x+\sum_{m=1}^\infty\fr{x^m}{m m!}
\label{6.18}
\eeq
is the exponential integral and $\Gg\approx 0.57721566$ is the Euler constant.

\begin{figure}[t]
\centerline{
\scalebox{0.7}{\includegraphics{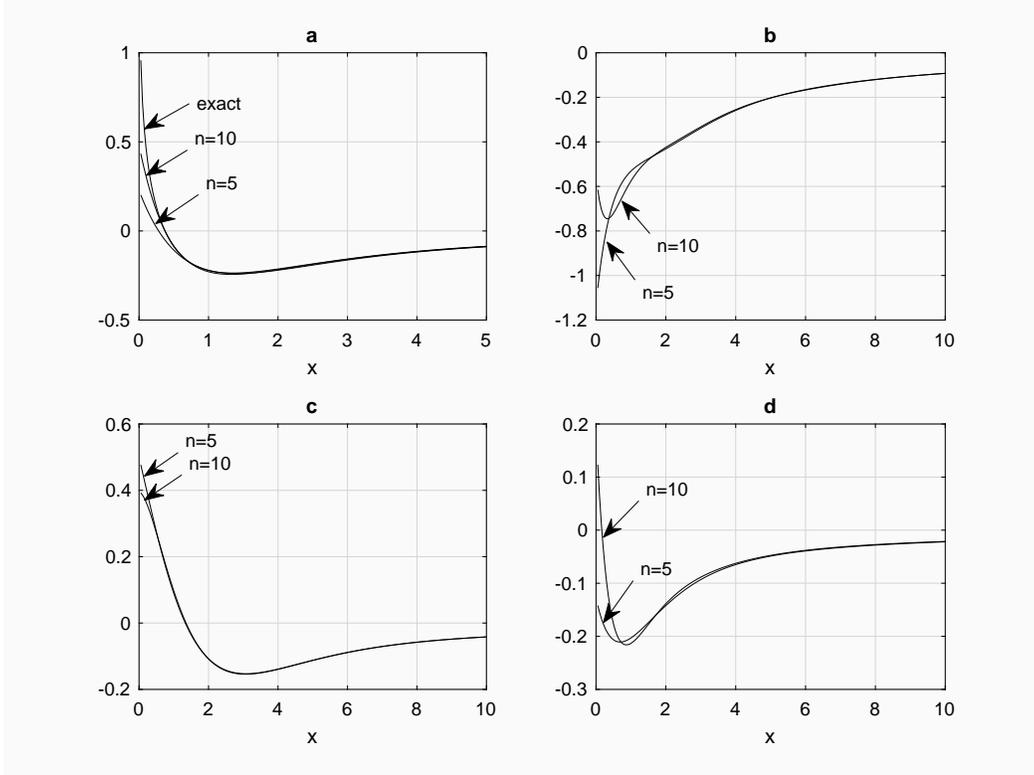}}}
\caption{The integral $I^\Ga[f](x)$ for $\Ga=-1/2$ when the number of zeros is $n=5$ and $n=10$.
(a): $f(x)=\sqrt{x}$ (in this case the integral is expressed through the exponential integral ${\rm Ei}(x)$). (b): $f(x)=e^{x/2}$.  (c): $f(x)=x\sqrt{x}$. (d): $f(x)=\sqrt{x}/(x+1)$.
}
\label{fig5}
\end{figure}

\setcounter{equation}{0}
  
\section{Conclusions}

We showed that the Hilbert transforms of the weighted Hermite polynomials $\exp(-x^2/2) H_n(x)$, the functions $G_n(x)$,
form a complete orthogonal system of functions in the space $L_2(-\infty,\infty)$. We 
discovered that the Hilbert transforms in a semi-axis of the weighted Laguerre polynomials $e^{-\Gn/2}\Gn^{-1/2}L_n^{-1/2}(x)$
and  $e^{-\Gn/2}\Gn^{1/2}L_n^{-1/2}(x)$
up to constant factors equal to the functions $\Gx^{-1/2}G_{2m}(\sqrt{\Gx})$ and  $G_{2m+1}(\sqrt{\Gx})$, respectively.
The system of functions $G_{2n}(x)$ and $G_{2n+1}$ may be employed for solving singular integral equations or their systems of the form
$S[\Gc(t)](x)+K[\Gc(t)](x)=f(x)$, $0<x<\infty$, in the class of integrable functions unbounded and bounded at the point $x=0$,
respectively. Here, $S$ is a singular operator with the Cauchy kernel, and $K$ is a regular operator. 
The method ultimately  reduces the integral equations to systems of linear algebraic equations of the second kind. If these systems are regular or at least quasiregular, then
 they can be solved numerically by the reduction method. The sufficiency of this scheme needs to be verified by means of numerical tests.
If $K=0$, then the solution is exact, and its representation is free of singular integrals. This method might be also employed 
for vector Riemann-Hilbert problems when the Wiener-Hopf factors are not available, and the associated system of    
integral equations has the structure $S[\Gc(t)](x)+K[\Gc(t)](x)=f(x)$, $0<x<\infty$. 

By employing an integral representation of the Jacobi function of the second kind $Q_n^{(\Ga,\Gb)}(x)$,
expressing it in terms of the hypergeometric Gauss function, and passing to the limit $n\to\infty$
in the representations for $n^{-\Ga}Q_n^{(\Ga,\Gb)}(1-\fr12z^2/n^2)$
we obtained the semi-infinite  Hilbert transforms of the Bessel function $J_\Ga(\Gl\sqrt{t})$ in terms of the functions $J_{\pm \Ga}(\Gl\sqrt{x})$ and $I_{\pm\Ga}(\Gl\sqrt{x})$
 in the intervals $0<x<\infty$ and $-\infty<x<0$, respectively. Here, $\Gl$ is a positive parameter.
We applied this result 
to derive a closed-form solution  to a model problem of
contact mechanics. The solution is free of singular integrals, and  the associated Riemann-Hilbert problem,
as the standard way of dealing with  such problems, was bypassed.

One of the most frequently applied methods for singular integral equations with the Cauchy kernel is the collocation method.
To employ it for singular integral equations in a semi-axis, one needs to make  the optimal choice of the collocation points
and have at their disposal an efficient procedure for the Cauchy integral in a semi-axis. To find such a formula, we proposed to use the Hilbert
transform of the weighted Laguerre polynomial $x^\Ga e^{-x}L_n^\Ga(x)$ derived in the paper, the Gauss quadrature 
formula for the integral $\int_0^\infty x^\Ga e^{-x}f(x)dx$ exact for polynomials of degree not higher than $2n-1$, and 
the Christoffel-Daurboux formula for the Laguerre polynomials. The quadrature formula for the singular integral with the 
Cauchy kernel $1/(t-x)$
in a semi-axis with the density $t^\Ga e^{-t}f(t)$ is exact for any plynomial $f(t)$ of degree not higher than $n-1$, requires computing
the Laguerre polynomials $L_n^\Ga(x)$ and $L_{n-1}^\Ga(x_m)$ and the confluent hypergeometric function $\GF(-n-\Ga,1-\Ga;t)$
at the points $t=x$ and $t=x_m$, where $x_m$ ($m=1,2,\ldots,n$) are the $n$ zeros of the polynomial  $L_n^\Ga(x)$.
Numerical tests proved efficiency of the quadrature formula.

\vspace{.2in}


\bibliographystyle{amsplain}

\begin{thebibliography}{10}

\bibitem{ant1}
Y.A. Antipov, \textit{Weight functions of a crack in a two-dimensional micropolar solid,}  Quart. J. Mech. Appl. Math., \textbf{65}, (2012), 239--271.

\bibitem{ant2} Y.A. Antipov and A.V. Smirnov, 
\textit{Subsonic propagation of a crack parallel to the boundary of a half-plane}, Math. Mech. Solids, \textbf{ 18} (2013), 153--167.

\bibitem{ant3} Y.A. Antipov,  \textit{Subsonic frictional cavitating penetration of a thin rigid body into an elastic medium},  Quart. J. Mech. Appl. Math.,  \textbf{71} (2018). 

\bibitem{bat1} H. Bateman,  \textit{Higher Transcendental Functions}, vol. 1, Bateman Manuscript Project, McGraw-Hill, New York, 1953.

\bibitem{bat2} H. Bateman,  \textit{ Higher Transcendental Functions}, vol. 2, Bateman Manuscript Project, McGraw-Hill, New York, 1953.

\bibitem{bat3} H. Bateman,  \textit{Tables of Integral Transforms}, vol. 1, Bateman Manuscript Project, McGraw-Hill, New York, 1954.

\bibitem{ell} J. Elliott,  \textit{On a class of integral equations,}  Proc. Amer. Math. Soc., 3  (1952),  566--572.

\bibitem{gak} F.D. Gakhov,   \textit{ Boundary Value Problems},  Pergamon Press,  Oxford, 1966.

\bibitem{gal} L.A. Galin,  \textit{Contact Problems in the Theory of Elasticity and Viscoelasticity} (Russian), Nauka, Moscow, 1980.

\bibitem{gra} I.S. Gradshte\u in and I.M.  Ryzhik, \textit{Table of Integrals, Series and Products,} Academic,  Oxford, Boston, 2007.

\bibitem{kop} W. Koppelman and J.D. Pincus,  \textit{Spectral representations for finite Hilbert transformations},  Math. Zeitschr., \textbf{71} (1959), 399--407.

\bibitem{kor} A.A. Korne\u i\v cuk,  \textit{Quadrature formulae for singular integrals} (Russian), {\v Zb. Vy\v cisl. Mat. i Mat. Fiz.,} \textbf{4}, no. 4, suppl. (1964), 64--74. 

\bibitem{lig}   M.J. Lighthill,  \textit{Introduction to Fourier Analysis and Generalized Functions,} Cambridge University Press,  Cambridge,  1964.                                                                           

\bibitem{mkh1}  S.M. Mkhitaryan, \textit{The spectral decompositions of integral operators which are analogous to the finite Hilbert transform},  Mat. Issled.,
 \textbf{4},  issue 1 (11) (1969), 98--109. 

\bibitem{mkh2}  S.M. Mkhitaryan, M.S. Mkrtchyan and E.G. Kanetsyan, 
\textit{On a method for solving Prandtl's integro-differential equation applied to problems of continuum mechanics using polynomial approximations},
ZAMM Z. Angew. Math. Mech.,
\textbf {97} (2017), 639--654. 

\bibitem{mus1} N.I. Muskhelishvili,   \textit{Singular Integral Equations},   P. Noordhoff,  Groningen, 1953.

\bibitem{mus2} N.I. Muskhelishvili,  \textit{Some Basic Problems of the Mathematical Theory of Elasticity},
 P. Noordhoff, Groningen, 1963.

\bibitem{pop} G.Ia. Popov,   \textit{Plane contact problem of the theory of elasticity with bonding or frictional forces}, J. Appl. Math.  Mech., \textbf{30} (1967), 653--667.

\bibitem{rei}  E. Reissner,   \textit{Boundary value problems in aerodynamics of lifting surfaces in non-uniform motion},  Bull. Amer. Math. Soc.,  {\textbf 55} (1949), 825--850.

\bibitem{son} H. S\"onngen,  \textit{Zur Theorie der endlichen Hilbert-Transformation,}   Math. Zeitschr.,  \textbf{ 60}  (1954), 31--51.

\bibitem{sta} V.J.E. Stark,   \textit{A generalized quadrature formula for Cauchy integrals},  AIAA Journal,  \textbf{9} (1971), 1854--1855.

\bibitem{sze} G. Szeg\"o,  \textit{Orthogonal Polynomials}, Amer. Math. Soc. Colloq. Publ., vol. 23, Providence, 1975. 

\bibitem{tit} E.C. Titchmarsh,  \textit{Introduction to the Theory of Fourier Integrals}, Clarendon Press, Oxford, 1948.

\bibitem{tri1} F.G. Tricomi,  \textit{On the finite Hilbert transformation},  Quart. J. Math.,  \textbf{2} (1951), 199--211.

\bibitem{tri2} F.G Tricomi,  \textit{ Integral Equations},  Interscience Publishers,  New York, 1957.




\end{thebibliography}

\end{document}